\documentclass[11pt,a4paper]{article}
\usepackage[T1]{fontenc}
\usepackage[utf8]{inputenc}
\usepackage{lmodern,amsmath,amssymb,amsthm,mathtools}
\usepackage[margin=27mm]{geometry}
\usepackage{microtype,booktabs,array,enumitem,xcolor}
\usepackage[hypertexnames=false,colorlinks=true,linkcolor=blue!40!black,citecolor=blue!40!black,urlcolor=blue!40!black]{hyperref}
\setlist{itemsep=3pt,topsep=5pt}
\allowdisplaybreaks
\numberwithin{equation}{section}
\newtheorem{theorem}{Theorem}[section]
\newtheorem{proposition}[theorem]{Proposition}
\newtheorem{lemma}[theorem]{Lemma}
\newtheorem{corollary}[theorem]{Corollary}
\theoremstyle{definition}
\newtheorem{definition}[theorem]{Definition}
\theoremstyle{remark}
\newtheorem{remark}[theorem]{Remark}
\newcommand{\R}{\mathbb R}
\newcommand{\RN}{\mathbb R^N}
\newcommand{\dx}{\,\mathrm dx}
\newcommand{\ii}{\int_{\RN}}
\newcommand{\A}{\mathcal A}
\newcommand{\E}{\mathcal E}
\newcommand{\Neh}{\mathcal N}
\newcommand{\Hr}{H^1_{\mathrm{rad}}(\RN)}
\newcommand{\X}{\mathcal X}
\newcommand{\nod}{\mathrm{nod}}
\newcommand{\dist}{\operatorname{dist}}
\newcommand{\diag}{\operatorname{diag}}

\newcommand{\one}{\boldsymbol 1}
\title{Radial nodal action ground states with two prescribed masses\newline for coupled Schr\"odinger equations}
\author{Yumeng Chou, Huabo Zhang}
\date{13 November 2026}
\hypersetup{pdftitle={Radial nodal action ground states with two prescribed masses for coupled Schrodinger equations},pdfauthor={Yumeng Chou, Huabo Zhang}}
\begin{document}
\maketitle
\begin{abstract}
We study normalized nodal solutions of a cooperative system of nonlinear Schrödinger equations with two independently prescribed masses. For sufficiently small positive coupling, we obtain a radial nodal action ground state for every mass pair in a fixed compact subset of the positive quadrant. We first minimize the action on a radial nodal Nehari set at fixed frequencies. We then select both frequencies through an interior minimization of the action level. In this way, we adapt the scalar action approach of De Coster--Dovetta--Galant--Serra (\emph{Ann. Inst. H. Poincaré Anal. Non Linéaire}, 2026) to construct solutions whose components both change sign and whose action is minimal among radial nodal solutions at the selected frequencies.
\end{abstract}

\section{Introduction and main results}\label{sec:intro}
We study the following system of coupled Schr\"odinger equations with prescribed masses:
\begin{equation}\label{system}
\begin{cases}
-\Delta u+\lambda_1u=\mu_1|u|^{p-2}u+\nu\alpha|u|^{\alpha-2}u|v|^\beta,
&x\in\RN,\\
-\Delta v+\lambda_2v=\mu_2|v|^{q-2}v+\nu\beta|u|^\alpha|v|^{\beta-2}v,
&x\in\RN,\\
u,v\in H^1(\RN),\qquad \displaystyle\ii u^2\dx=a,\quad
\displaystyle\ii v^2\dx=b,
\end{cases}
\end{equation}
in which $N\geq3$, $u,v$ are real-valued functions, $\mu_1,\mu_2>0$ and $\nu\geq0$. The masses $a,b>0$ are prescribed, while the frequencies $\lambda_1,\lambda_2$ are unknown. The power exponents $p,q>2$, and the coupling exponents $\alpha,\beta>1$. For $1\leq t<\infty$, we write $\|\cdot\|_t$ for the norm in $L^t(\RN)$. We denote the Sobolev critical exponent by $2^*=\frac{2N}{N-2}$, the mass critical exponent by $p_c=2+\frac{4}{N}$, and the interaction degree by $r=\alpha+\beta$. The energy and action are
\begin{align}
\E_\nu(u,v)&=\frac12\ii(|\nabla u|^2+|\nabla v|^2)\dx
-\frac{\mu_1}{p}\ii|u|^p\dx-\frac{\mu_2}{q}\ii|v|^q\dx
-\nu\ii|u|^\alpha|v|^\beta\dx,\label{energy}\\*
\A_{\lambda,\nu}(u,v)&=\E_\nu(u,v)+\frac{\lambda_1}{2}\|u\|_2^2+
\frac{\lambda_2}{2}\|v\|_2^2,\qquad \lambda=(\lambda_1,\lambda_2).
\label{action}
\end{align}
We call a pair nodal $(u,v)$ if both components change sign.\textbf{We look for radial nodal solutions of Equations~\eqref{system} that minimize the action among radial nodal solutions at the same frequencies.}  We first treat the frequencies as parameters and then select them so that the two mass constraints hold.

For scalar Schr\"odinger equations, Dovetta-Serra-Tilli~\cite{DST2023} established a Legendre--Fenchel relation between the action and energy ground-state levels and studied the dependence of the masses of action ground states on the frequency. De Coster-Dovetta-Galant-Serra~\cite{DCDGS} used the action level to study normalized nodal solutions on bounded domains. They characterized the masses of nodal action ground states and determined regimes in which least energy normalized nodal solutions are nodal action ground states. We consider a system on the whole space, with two frequencies and two prescribed masses.

For systems at fixed frequencies, Liu-Liu-Wang~\cite{LLWnodal} developed a variational framework for mixed nodal states and obtained infinitely many solutions in the repulsive-coupling regime, with every component changing sign. Li-Wang~\cite{LiWang} constructed infinitely many solutions with prescribed componentwise nodal numbers for an $m$-component cubic system on radial domains under a repulsive-coupling condition. For a two-component system in $\mathbb R^3$, Zhang~\cite{Zhang} obtained arbitrarily many radial sign-changing solutions. For sufficiently small positive coupling, the author also obtained a least energy radial sign-changing solution. The powers in that model have degrees in $[4,6)$, whereas the powers in our normalization result have degrees below $\frac{10}{3}$.

On bounded smooth domains in dimensions $N=2,3$, Amadori~\cite{Amadori} studied a two-component system with a common Sobolev subcritical exponent $2q$. The author characterized least energy nodal and seminodal solutions as minima on Nehari-type sets of codimensions four and three, computed their Morse indices and nodal domains, and proved nonradiality on balls for sufficiently small coupling. Bieganowski~\cite{Bieganowski} studied systems with sign-changing nonlinearities by using the  concentration-compactness and Lusternik--Schnirelmann theory.

For prescribed-mass solutions of coupled systems, Gou-Jeanjean~\cite{GJ2016} obtained constrained global minimizers in the fully $L^2$-subcritical regime for arbitrary positive masses and coupling. Bartsch-Jeanjean-Soave~\cite{BJS2016} studied the cubic two-component system in $\mathbb R^3$ and obtained positive normalized solutions for coupling ranges depending on the prescribed masses. Bartsch-Soave~\cite{BS2017} used a natural constraint for the competing cubic system in $\mathbb R^3$ and studied phase separation as the coupling tends to negative infinity. Bartsch-Jeanjean~\cite{BJ2018} treated more general homogeneous nonlinearities, including $L^2$-supercritical powers, by constrained linking.

In mass-mixed regimes, Gou-Jeanjean~\cite{GJ2018} obtained two positive normalized solutions for sufficiently small attractive coupling under the corresponding exponent assumptions. One solution is a local minimizer; the other is found by a constrained mountain-pass or linking argument. Writing $p_1,p_2$ for their power exponents and $r_1,r_2$ for their interaction exponents, we recall that their first regime requires
\[
2<p_1,p_2<2+\frac{4}{N}<r_1+r_2<2^*.
\]
In dimension $N=3$, this includes the arrangement $p,q<p_c<r<2^*$ in our main theorem. In dimensions $N\geq5$, their result for the second solution requires an additional restriction on the exponents.

Normalized nodal solutions of systems have also been constructed on bounded domains. Liu-Song-Wu-Zou~\cite{LiuSongWuZou} considered an $m$-component cubic Gross--Pitaevskii system in dimensions $N=3,4$, with nonzero self-interaction coefficients and symmetric nonzero inter-component coefficients of arbitrary signs. For every positive integer $j$, they obtained at least $j$ sign-changing normalized solutions when their prescribed mass parameters $c_1,\ldots,c_m>0$ satisfy the following condition, where $\widetilde c_j>0$ is their threshold:
\[
\frac{\max_{1\leq i\leq m}c_i^2}{\min_{1\leq i\leq m}c_i}<\widetilde c_j.
\]
They proved a corresponding multiplicity result for seminodal solutions, using vector linking and partial vector linking, and studied bifurcation as the mass vector tends to zero. Unlike this work, we work on the whole space and seek a radial nodal solution characterized by an action minimum.

For critical nonlinearities, Zhang-Zhong-Zhou~\cite{ZZZ} considered $p=q=2^*$ in dimensions $N=3,4$ with $\alpha+\beta<2+\frac{4}{N}$. For sufficiently small coupling and arbitrary fixed positive masses, they obtained a positive radial local minimizer and a second positive radial solution of mountain-pass type. Guo-He-Shuai-Zhong~\cite{GHShZ} studied
\[
2<p,q<2+\frac4N,\qquad \alpha,\beta>1,\qquad \alpha+\beta=2^*,\qquad N\geq3.
\]
For every fixed positive mass pair and sufficiently small coupling, they constructed a positive normalized local minimizer and a positive normalized mountain-pass solution. Their coupling bound depends on upper bounds for the masses.

We first fix the frequencies and minimize the action over the radial nodal Nehari set, imposing one constraint for each of the four sign parts. Under theassumption of~\eqref{H} and for sufficiently small nonnegative coupling, we prove that the infimum is attained and that every minimizer is a full weak solution. The coupling bound is uniform for frequencies in any compact rectangle in the positive quadrant. In particular, each minimizer is an action ground state among radial nodal solutions at those frequencies. We assume $\alpha,\beta\geq2$ and $\alpha+\beta\leq2^*$, so we have $N=3,4$. When $N=4$ we must have $\alpha=\beta=2$ and $r=2^*$.

To prescribe the masses, we assume in addition that $2<p,q<2+\frac{4}{N}$. The powers are then mass-subcritical, while the interaction is mass-supercritical and may be Sobolev-critical. Given any compact rectangle of positive mass pairs, we select the two frequencies by minimizing a function of the action level. For sufficiently small positive coupling, with a bound depending on the rectangle, we obtain a radial nodal action ground state for each prescribed mass pair. The action-minimizing property holds in the radial class at the selected frequencies.

We impose the following assumptions for the fixed-frequency construction:
\begin{equation}\tag{H}\label{H}
N\geq3,\qquad 2<p,q<2^*,\qquad \alpha,\beta\geq2,
\qquad r\leq2^*,\qquad \mu_1,\mu_2>0.
\end{equation}
For the normalization result, we additionally assume
\begin{equation}\tag{M}\label{M}
2<p,q<p_c.
\end{equation}
In Sections~\ref{sec:bounds}--\ref{sec:normalization} we work Under theassumption of~\eqref{H}, adding~\eqref{M} in Section~\ref{sec:normalization}. In Section~\ref{sec:criticalpure}, we use different assumptions on the exponents.

Write $\Hr$ for the radial subspace of $H^1(\RN)$. For $w\in H^1(\RN)$, let $w^+=\max\{w,0\}$ and $w^-=\max\{-w,0\}$; for $\ell>0$, set $\|w\|_{1,\ell}^2=\|\nabla w\|_2^2+\ell\|w\|_2^2$.

For positive frequencies and $\sigma,\tau\in\{+,-\}$, define the four constraint functions
\begin{align}
G_{u,\sigma}(u,v)&=\|u^\sigma\|_{1,\lambda_1}^2-
\mu_1\|u^\sigma\|_p^p-
\nu\alpha\ii|u^\sigma|^\alpha|v|^\beta\dx,\label{Gu}\\
G_{v,\tau}(u,v)&=\|v^\tau\|_{1,\lambda_2}^2-
\mu_2\|v^\tau\|_q^q-
\nu\beta\ii|u|^\alpha|v^\tau|^\beta\dx.\label{Gv}
\end{align}
We denote the radial four-sign nodal Nehari set by $\Neh_{\lambda,\nu}$ and its action level by $m_\nu(\lambda)$:
\begin{align}
\Neh_{\lambda,\nu}
&=\{(u,v)\in\Hr\times\Hr:u^\pm,v^\pm\ne0,\quad
G_{u,+}=G_{u,-}=G_{v,+}=G_{v,-}=0\},\label{Nnod}\\
m_\nu(\lambda)&=\inf_{\Neh_{\lambda,\nu}}
\A_{\lambda,\nu}.\label{level}
\end{align}
A full weak solution satisfies both equations when tested against every function in $H^1(\RN)$.

\begin{definition}\label{def:ground}
A radial nodal action ground state at frequencies $\lambda$ is a radial nodal solution of the first two equations in Equations~\eqref{system} with smallest action among all radial nodal solutions at the same frequencies. When it also has masses $(a,b)$, we call it a normalized radial nodal action ground state.
\end{definition}
\begin{theorem}\label{thm:fixed}
Assume~\eqref{H}. For every compact rectangle $K\subset(0,\infty)^2$, there is $\nu_K>0$ such that, for $0\leq\nu<\nu_K$ and $\lambda\in K$, the level $m_\nu(\lambda)$ is positive and attained. Every minimizer is a full weak solution of the first two equations in Equations~\eqref{system}. In particular,
\begin{equation}\label{least-solutions}
m_\nu(\lambda)=\min\{\A_{\lambda,\nu}(u,v):
(u,v)\in\Hr\times\Hr\text{ is a nodal solution at }\lambda\}.
\end{equation}
For fixed $0\leq\nu<\nu_K$, the set of these minimizers is compact in $H^1(\RN)^2$, also along convergent frequency sequences in $K$. The level is locally Lipschitz in $\lambda$ on the interior of $K$.
\end{theorem}

For a scalar power $s\in(2,2^*)$, coefficient $\mu>0$, and frequency $\ell>0$, denote by $d_{s,\mu}(\ell)$ the infimum of
\[
I_{s,\mu,\ell}(w)=\frac12\|w\|_{1,\ell}^2-\frac\mu s\|w\|_s^s
\]
over radial $w$ with $w^\pm\ne0$ and
$\|w^\pm\|_{1,\ell}^2=\mu\|w^\pm\|_s^s$.
We prove attainment in Section~\ref{sec:parameter}. Set
\begin{equation}\label{sigmakappa}
\sigma_1=\frac{2}{p-2}-\frac N2,\quad
\sigma_2=\frac{2}{q-2}-\frac N2,\quad
d_1=d_{p,\mu_1}(1),\quad d_2=d_{q,\mu_2}(1).
\end{equation}
Under~\eqref{M}, both $\sigma_1$ and $\sigma_2$ are positive. We write $\operatorname{int} K$ and $\partial K$ for the interior and boundary of $K$.

\begin{theorem}\label{thm:mass}
Assume~\eqref{H} and~\eqref{M}. Fix positive endpoints $0<a_-\leq a_+$, $0<b_-\leq b_+$ and the mass rectangle
$\mathcal M=[a_-,a_+]\times[b_-,b_+]$.
There exist a compact rectangle $K\subset(0,\infty)^2$ and $\nu_{\mathcal M}>0$ such that, for every $(a,b)\in\mathcal M$ and $0<\nu<\nu_{\mathcal M}$, there are frequencies $\lambda_\nu\in\operatorname{int} K$ and a normalized radial nodal action ground state $(u_\nu,v_\nu)$ satisfying Equations~\eqref{system}. The masses satisfy
\begin{equation}\label{mass-exact}
\|u_\nu\|_2^2=a,\qquad \|v_\nu\|_2^2=b.
\end{equation}
The frequencies can be selected to minimize
\begin{equation}\label{Fdef}
F_{\nu,a,b}(\lambda)=m_\nu(\lambda)-\frac12(a\lambda_1+b\lambda_2)
\quad\text{on }K.
\end{equation}
At any such interior minimizer every action minimizer has the masses in~\eqref{mass-exact}, and
\begin{equation}\label{gradient}
m_\nu\text{ is differentiable at }\lambda_\nu,
\qquad 2\nabla m_\nu(\lambda_\nu)=(a,b).
\end{equation}
As $\nu\downarrow0$, uniformly for $(a,b)\in\mathcal M$,
\begin{equation}\label{lambda-limit}
\lambda_{\nu,j}\longrightarrow\lambda_j^0,
\quad \lambda_1^0=\left(\frac{a}{2(\sigma_1+1)d_1}\right)^{\frac{1}{\sigma_1}},
\quad \lambda_2^0=\left(\frac{b}{2(\sigma_2+1)d_2}\right)^{\frac{1}{\sigma_2}}.
\end{equation}
Write $\lambda^0=(\lambda_1^0,\lambda_2^0)$ for this limiting frequency pair. For fixed masses, subsequences of the states converge strongly in $H^1(\RN)^2$ to products of scalar radial nodal action ground states at $\lambda^0$. Moreover,
\begin{equation}\label{negative-limit}
\E_\nu(u_\nu,v_\nu)\longrightarrow
-\frac{\sigma_1}{2(\sigma_1+1)}a\lambda_1^0
-\frac{\sigma_2}{2(\sigma_2+1)}b\lambda_2^0<0.
\end{equation}
\end{theorem}

\begin{remark}
In Section~\ref{sec:criticalpure}, we consider critical powers with a subcritical interaction under assumptions different from~\eqref{H}--\eqref{M}. Under the coefficient-balance condition, we reduce the system to a scalar equation using the proportionality relation $(u,v)=(\sqrt{\alpha}\,w,\sqrt{\beta}\,w)$, and the mass pairs therefore lie on the ray $(a,b)=(\alpha m,\beta m)$, where $m=\|w\|_2^2$. We minimize the action within this proportionality class. In the Sobolev-critical case, the interaction range is $2<r<p_c$, a proper subrange of $2<r<2^*$. The positive normalized solutions of~\cite{ZZZ} allow independent masses, whereas we obtain sign-changing states on a restricted mass ray with the stated action-minimizing property.

\end{remark}

We organize the paper as follows. In Section~\ref{sec:bounds}, we collect notation and establish uniform lower and upper bounds on the four nodal components. In Section~\ref{sec:deform}, we construct a bounded Palais--Smale sequence by a four-amplitude deformation. In Section~\ref{sec:compact}, we prove compactness and attainment of the fixed-frequency radial nodal action level. In Section~\ref{sec:parameter}, we study its dependence on the two frequencies and derive the mass identities at an interior minimum. In Section~\ref{sec:normalization}, we complete the proof of Theorem~\ref{thm:mass}. Finally, in Section~\ref{sec:criticalpure}, we prove the existence results for critical powers in the proportionality class.

\section{Preliminaries and bounds on the nodal Nehari set}\label{sec:bounds}
We set $\X=\Hr\times\Hr$ and denote its dual space by $\X^*$. We equip $\X$ with the product $H^1$ norm. Integrals are taken over $\RN$ unless otherwise indicated. With the sign convention introduced in Section~\ref{sec:intro}, we have $w=w^+-w^-$ and $w^+w^-=0$ almost everywhere.

The action is continuously Fr\'echet differentiable on $H^1(\RN)^2$ in the exponent ranges considered below. We denote its derivative by a prime and write
\begin{equation}\label{mixedintegral}
C(u,v)=\ii|u|^\alpha|v|^\beta\dx.
\end{equation}

Fix a compact positive frequency rectangle $K$. Constants denoted by $c,C>0$ may change from line to line; unless additional dependence is stated, they depend only on the fixed exponents, coefficients, $K$, and the chosen action bound. They are independent of the functions and of sufficiently small $\nu$. The symbols $O(t)$ and $o(1)$ denote, respectively, a quantity bounded by a fixed multiple of $|t|$ and a quantity tending to zero in the stated limit. For the four sign parts, we introduce the coefficients
\begin{align*}
A_\sigma&=\|u^\sigma\|_{1,\lambda_1}^2,&
B_\sigma&=\mu_1\|u^\sigma\|_p^p,&
T_{\sigma\tau}&=\nu\ii|u^\sigma|^\alpha|v^\tau|^\beta\dx,\\
D_\tau&=\|v^\tau\|_{1,\lambda_2}^2,&
E_\tau&=\mu_2\|v^\tau\|_q^q,&
\nu C(u,v)&=\sum_{\sigma,\tau}T_{\sigma\tau}.
\end{align*}
The nodal constraints can then be written as
\begin{equation}\label{coeffconstraints}
A_\sigma=B_\sigma+\alpha\sum_\tau T_{\sigma\tau},
\qquad D_\tau=E_\tau+\beta\sum_\sigma T_{\sigma\tau}.
\end{equation}
All sums over sign indices run over $\{+,-\}$. Set $Q=\sum_\sigma A_\sigma+\sum_\tau D_\tau=\|u\|_{1,\lambda_1}^2+\|v\|_{1,\lambda_2}^2$ and $k_t=\frac{1}{2}-\frac{1}{t}$ for $t>2$. On the nodal Nehari set,
\begin{equation}\label{coercive}
Q=\sum_\sigma B_\sigma+\sum_\tau E_\tau+\nu rC(u,v),\qquad
\A_{\lambda,\nu}=k_p\sum_\sigma B_\sigma+k_q\sum_\tau E_\tau+\nu rk_rC(u,v)
\geq k_*Q,
\end{equation}
where $k_*=\min\{k_p,k_q,k_r\}>0$.

\begin{lemma}\label{lem:trial}
For sufficiently small $\nu\geq0$, the set~\eqref{Nnod} is nonempty and $m_\nu(\lambda)\leq L_0$ uniformly on $K$, for some finite $L_0$.
\end{lemma}
\begin{proof}
Choose smooth compactly supported radial functions $f,g$ with both signs nonzero. At $\nu=0$, the required positive amplitudes are
\[
s_\sigma^0=\left(\frac{\|f^\sigma\|_{1,\lambda_1}^2}
{\mu_1\|f^\sigma\|_p^p}\right)^{\frac{1}{p-2}},\qquad
t_\tau^0=\left(\frac{\|g^\tau\|_{1,\lambda_2}^2}
{\mu_2\|g^\tau\|_q^q}\right)^{\frac{1}{q-2}}.
\]
The amplitudes are continuous and bounded above and away from zero for $\lambda\in K$. The Jacobian of the four constraints with respect to these amplitudes is diagonal at $\nu=0$, with entries
$-(p-2)s_\sigma^0\|f^\sigma\|_{1,\lambda_1}^2$ and
$-(q-2)t_\tau^0\|g^\tau\|_{1,\lambda_2}^2$. These entries are bounded away from zero uniformly on $K$. By the finite-dimensional implicit-function theorem, we find nearby positive amplitudes for all sufficiently small $\nu$, uniformly on $K$. The resulting trial actions are uniformly bounded.
\end{proof}

Fix $L>L_0+1$ and consider points of~\eqref{Nnod} with action at most $L$.
By~\eqref{coercive}, $Q\leq M:=\frac{L}{k_*}$. Using the embeddings of $H^1(\RN)$ into $L^r(\RN)$ and the lower bounds for the frequencies, we obtain
\begin{align}
\nu\alpha\ii|u^\sigma|^\alpha|v|^\beta\dx
&\leq\nu\alpha\|u^\sigma\|_r^\alpha\|v\|_r^\beta
\leq C\nu A_\sigma^{\frac{\alpha}{2}}Q^{\frac{\beta}{2}}\nonumber\\
&=C\nu A_\sigma A_\sigma^{\frac{\alpha-2}{2}}Q^{\frac{\beta}{2}}
\leq C\nu M^{\frac{r-2}{2}}A_\sigma.\label{absorb}
\end{align}
In the last inequality, we used $\alpha\geq2$. For the analogous inequality for a sign part of $v$, we use $\beta\geq2$.

\begin{lemma}\label{lem:nonvanish}
After decreasing the upper bound for $\nu$, there are constants $c,C>0$ such that every point in~\eqref{Nnod} of action at most $L$ satisfies
\begin{equation}\label{uniformcoeff}
c\leq A_\sigma,B_\sigma,D_\tau,E_\tau\leq C,
\qquad 0\leq T_{\sigma\tau}\leq C\nu.
\end{equation}
The constants are independent of the nodal sets of the functions.
\end{lemma}
\begin{proof}
Choose $\nu$ so that the coefficient in~\eqref{absorb} is at most $\frac{1}{2}$. From the constraint, we obtain
$B_\sigma\geq \frac{A_\sigma}{2}$, while by Sobolev embedding we have
$B_\sigma\leq C A_\sigma^{\frac{p}{2}}$. Combining these inequalities and dividing by $A_\sigma>0$, we find
\[
\frac12\leq C A_\sigma^{\frac{p-2}{2}},\qquad
A_\sigma\geq(2C)^{-\frac{2}{p-2}}>0.
\]
From $B_\sigma\geq \frac{A_\sigma}{2}$, we obtain the lower bound for $B_\sigma$. We obtain lower bounds for $D_\tau$ and $E_\tau$ in the same way. We obtain the upper bounds from $Q\leq M$ and~\eqref{coeffconstraints}; by H\"older's inequality and Sobolev embedding, we also have $T_{\sigma\tau}\leq C\nu$.
\end{proof}

\section{Four amplitudes and a Palais--Smale sequence}\label{sec:deform}
A Palais--Smale sequence at level $c$ is a sequence whose action tends to $c$ and whose action derivative tends to zero in the dual norm. We will also establish boundedness in $\X$. For $z=(u,v)$ in~\eqref{Nnod} and an amplitude vector $\xi=(s_+,s_-,t_+,t_-)\in[0,\infty)^4$, put
\[
\gamma_z(s_+,s_-,t_+,t_-)
=(s_+u^+-s_-u^-,\ t_+v^+-t_-v^-)
\]
and $J_z=\A_{\lambda,\nu}\circ\gamma_z$. On $[0,\infty)^4$,
\begin{equation}\label{fiber}
J_z(\xi)=\frac12\sum_\sigma A_\sigma s_\sigma^2
+\frac12\sum_\tau D_\tau t_\tau^2
-\frac1p\sum_\sigma B_\sigma s_\sigma^p
-\frac1q\sum_\tau E_\tau t_\tau^q
-\sum_{\sigma,\tau}T_{\sigma\tau}s_\sigma^\alpha t_\tau^\beta.
\end{equation}
Write $\one=(1,1,1,1)$.

\begin{lemma}\label{lem:fiber}
Uniformly for the points of Lemma~\ref{lem:nonvanish}, and for sufficiently small $\nu$, $J_z$ has exactly one critical point in $(0,\infty)^4$. It is $\one$, and it is the unique global maximum on $[0,\infty)^4$. The Hessian at this point satisfies
\begin{equation}\label{Hnegative}
D^2J_z(\one)\leq-\eta I_4
\end{equation}
where $D^2J_z$ denotes the Hessian, $I_4$ is the $4\times4$ identity matrix, and $\eta>0$ is uniform. On $\mathcal B=[\frac{1}{2},\frac{3}{2}]^4$, there is a uniform constant $c_{\mathrm{gap}}>0$ such that
\begin{equation}\label{uniformquadraticgap}
J_z(\one)-J_z(\xi)\geq c_{\mathrm{gap}}|\xi-\one|^2
\qquad\text{for every }\xi\in\mathcal B.
\end{equation}
\end{lemma}
\begin{proof}
Since $p,q>2$ and the mixed terms have negative sign, $J_z(\xi)\to-\infty$ as $|\xi|\to\infty$ in the positive orthant, uniformly for the coefficients in~\eqref{uniformcoeff}. Thus a global maximum exists. At a maximum, every positive coordinate satisfies its partial critical-point equation, and in particular
\[
B_\sigma s_\sigma^{p-2}\leq A_\sigma,
\qquad E_\tau t_\tau^{q-2}\leq D_\tau.
\]
All active coordinates are therefore uniformly bounded. If $s_\sigma=0$, increasing it to $\delta>0$ changes the action by
\[
\frac12 A_\sigma\delta^2-\frac1p B_\sigma\delta^p
-\delta^\alpha\sum_\tau T_{\sigma\tau}t_\tau^\beta.
\]
This is positive for sufficiently small $\delta$: if $\alpha>2$, the first term dominates, and if $\alpha=2$, its coefficient stays positive by~\eqref{uniformcoeff} and smallness of $\nu$. We treat a zero $t_\tau$ in the same way. Every maximum is consequently interior.

At any interior critical point,
\begin{equation}\label{fibercritical}
A_\sigma=B_\sigma s_\sigma^{p-2}
+\alpha s_\sigma^{\alpha-2}\sum_\tau T_{\sigma\tau}t_\tau^\beta,
\end{equation}
and analogously for $t_\tau$. Using the upper bounds obtained above, $\alpha,\beta\geq2$, and~\eqref{uniformcoeff}, we see that the last term is $O(\nu)$, including when a coordinate is close to zero. Combining this bound with~\eqref{coeffconstraints}, we find
$s_\sigma^{p-2}=1+O(\nu)$ and $t_\tau^{q-2}=1+O(\nu)$.
Thus all interior critical points lie in a fixed small box about $\one$ once $\nu$ is small.

At $\one$, the diagonal entries of the Hessian are
\begin{align}
\partial_{s_\sigma}^2J_z(\one)&=-(p-2)B_\sigma
-\alpha(\alpha-2)\sum_\tau T_{\sigma\tau},\label{Hdiag}\\
\partial_{t_\tau}^2J_z(\one)&=-(q-2)E_\tau
-\beta(\beta-2)\sum_\sigma T_{\sigma\tau},\nonumber
\end{align}
the mixed entries are
\begin{equation}\label{Hcross}
\partial_{s_\sigma}\partial_{t_\tau}J_z(\one)=-\alpha\beta T_{\sigma\tau},
\end{equation}
and the remaining off-diagonal entries vanish. For any vector $h\in\R^4$, the diagonal part is at most
$-c\min\{p-2,q-2\}|h|^2$, while the mixed part has absolute value at most $C\nu|h|^2$, by $2|h_i h_j|\leq h_i^2+h_j^2$. After decreasing $\nu$, we obtain~\eqref{Hnegative}. The same uniform negative definiteness holds throughout a sufficiently small box about $\one$: first choose the box using $c\leq B_\sigma,E_\tau\leq C$ and $p,q>2$, and then decrease $\nu$. By strict concavity on this box, we obtain uniqueness of the interior critical point. From the constraints, we have $\nabla J_z(\one)=0$, so this point is the unique global maximum and we obtain~\eqref{Hnegative}.

To prove the uniform gap~\eqref{uniformquadraticgap},
for $s>2$ put
\[
 f_s(x)=\left(\frac12-\frac1s\right)-\frac{x^2}{2}+\frac{x^s}{s}
 \qquad (x\geq0).
\]
The function $x\mapsto \frac{x^2}{2}-\frac{x^s}{s}$ has its unique positive maximum at $x=1$.
Consequently $f_s(x)>0$ for $x\ne1$, while
$f_s(1)=f_s'(1)=0$ and $f_s''(1)=s-2>0$.
In particular, with the value at $x=1$ defined by continuity,
\begin{equation}\label{fp-gap-constant}
 a_0:=\min\left\{
 \min_{\frac{1}{2}\leq x\leq\frac{3}{2}}\frac{f_p(x)}{(x-1)^2},
 \min_{\frac{1}{2}\leq x\leq\frac{3}{2}}\frac{f_q(x)}{(x-1)^2}
 \right\}>0.
\end{equation}
For each mixed pair define
\[
 g(s,t)=\left(\frac r2-1\right)-\frac\alpha2s^2
          -\frac\beta2t^2+s^\alpha t^\beta.
\]
We have $g(1,1)=0$ and $\nabla g(1,1)=0$. Since $g$ is $C^2$ on
$[\frac{1}{2},\frac{3}{2}]^2$, there is a number $M_g>0$, depending only on
$\alpha,\beta$, such that
\begin{equation}\label{g-quadratic}
 |g(s,t)|\leq \frac{M_g}{2}\big((s-1)^2+(t-1)^2\big)
 \qquad(\frac{1}{2}\leq s,t\leq\frac{3}{2}).
\end{equation}
Using~\eqref{coeffconstraints}, the difference of the fibre values is
\begin{align}\label{fiber-explicit-gap}
 J_z(\one)-J_z(\xi)
 ={}&\sum_\sigma B_\sigma f_p(s_\sigma)
     +\sum_\tau E_\tau f_q(t_\tau)
     +\sum_{\sigma,\tau}T_{\sigma\tau}g(s_\sigma,t_\tau).
\end{align}
By~\eqref{uniformcoeff}, the first two sums are at least
$c a_0|\xi-\one|^2$. Also $\sum_{\sigma,\tau}T_{\sigma\tau}\leq C_T\nu$
for a constant $C_T$ independent of $z$ and $\lambda\in K$.
By~\eqref{g-quadratic}, we bound the absolute value of the last sum
by $M_g C_T\nu|\xi-\one|^\frac{2}{2}$. We decrease the allowed coupling so that
$M_g C_T\nu\leq c a_0$. With
\[
 c_{\mathrm{gap}}=\frac{c a_0}{2}>0,
\]
we obtain~\eqref{uniformquadraticgap}. All constants are uniform for bounded-action points in the four-sign Nehari set and for frequencies in $K$.
\end{proof}

The constraint maps~\eqref{Gu}--\eqref{Gv} are continuous. In the deformation argument, we work in the ambient radial space $\X$ without differentiating the sign truncation maps.

\begin{lemma}\label{lem:PSsequence}
Fix a compact positive frequency rectangle $K$ and the sublevel bound $L$
from Section~\ref{sec:bounds}. After decreasing the permitted coupling,
let $\varepsilon_n\downarrow0$ and choose
$z_n=(u_n,v_n)\in\Neh_{\lambda,\nu}$ such that
\begin{equation}\label{chosen-min-seq}
 m_\nu(\lambda)\leq\A_{\lambda,\nu}(z_n)
 \leq m_\nu(\lambda)+\varepsilon_n^2.
\end{equation}
There are points $y_n\in\X$ and a constant $C>0$, uniform for
$\lambda\in K$ and small $\nu\geq0$, for which
\begin{equation}\label{PSsequence}
 \begin{split}
 |\A_{\lambda,\nu}(y_n)-m_\nu(\lambda)|&\leq2\varepsilon_n^2,\\
 \|(\A_{\lambda,\nu}|_\X)'(y_n)\|_{\X^*}&\leq8\varepsilon_n,\\
 \|y_n-z_n\|_\X&\leq C\varepsilon_n
 \end{split}
\end{equation}
for all sufficiently large $n$. In particular, $y_n$ is a bounded
Palais--Smale sequence at $m_\nu(\lambda)$, and its four sign parts
remain uniformly nonzero in the relevant $L^p$ or $L^q$ norms.
Every attained minimizer of $m_\nu(\lambda)$ is a critical point of
$\A_{\lambda,\nu}|_\X$.
\end{lemma}
\begin{proof}
Write $\A=\A_{\lambda,\nu}|_\X$, $m=m_\nu(\lambda)$, and
$\gamma_n=\gamma_{z_n}|_{\mathcal B}$.
By Lemmas~\ref{lem:trial} and~\ref{lem:nonvanish}, the points $z_n$
and the compact sets $\Gamma_n:=\gamma_n(\mathcal B)$ are uniformly
bounded in $\X$. The quadratic estimate~\eqref{uniformquadraticgap}
gives, for large $n$,
\begin{equation}\label{PS-boundary}
 \max_{\xi\in\partial\mathcal B}\A(\gamma_n(\xi))
 \leq m+\varepsilon_n^2-\frac{c_{\mathrm{gap}}}{4}
 <m-2\varepsilon_n^2.
\end{equation}
Here and below the norm on $\X$ is the standard product $H^1$ norm;
the corresponding dual norm is denoted by $\|\cdot\|_{\X^*}$.

\emph{1. A fixed open set preserving all four signs.}
From the lower bounds in~\eqref{uniformcoeff}, we find that
$\|u_n^\pm\|_p$ and $\|v_n^\pm\|_q$ are uniformly bounded below.
Since each amplitude on $\mathcal B$ is at least $\frac{1}{2}$, we can fix
$d>0$, independent of $n$, $\lambda$ and small $\nu$, such that
\begin{equation}\label{O-four}
 \|U^\pm\|_p\geq4d,\qquad\|V^\pm\|_q\geq4d
 \quad\text{for }(U,V)\in\Gamma_n.
\end{equation}
Define
\[
 \mathcal O=\{(U,V)\in\X:\|U^\pm\|_p>d,
                                    \ \|V^\pm\|_q>d\}.
\]
The truncation inequalities
$\|U^\pm-\widetilde U^\pm\|_p\leq\|U-\widetilde U\|_p$
and their $v$ counterparts, together with Sobolev embedding, imply
that there is $d_0>0$, independent of $n$, such that, with $\operatorname{dist}$ denoting distance in $\X$,
\begin{equation}\label{O-tube}
 \{y\in\X:\dist(y,\Gamma_n)\leq d_0\}\subset\mathcal O.
\end{equation}
In particular $G=(G_{u,+},G_{u,-},G_{v,+},G_{v,-})$ is a continuous
map from $\mathcal O$ into $\R^4$. This step uses continuity of
$U\mapsto U^\pm$ in $H^1$ and in Lebesgue spaces, not differentiability
of the truncation map.

\emph{2. The degree of the undeformed four-amplitude family.}
For $\xi\in\mathcal B$ we have, by direct differentiation in the
four positive real amplitudes,
\begin{equation}\label{G-on-path}
 G(\gamma_n(\xi))=\diag(\xi)\nabla J_{z_n}(\xi).
\end{equation}
Here $\operatorname{diag}(\xi)$ is the diagonal matrix with diagonal entries given by $\xi$. The right side has the unique zero $\one$ in $\mathcal B$.
The Jacobian there is $D^2J_{z_n}(\one)$, because
$\nabla J_{z_n}(\one)=0$. By~\eqref{Hnegative}, this symmetric matrix
has four strictly negative eigenvalues. Its determinant is positive;
therefore, writing $\deg$ for the Brouwer degree,
\begin{equation}\label{PS-degree}
 \deg(G\circ\gamma_n,\operatorname{int}\mathcal B,0)=+1.
\end{equation}
For any continuous $\widetilde\gamma:\mathcal B\to\mathcal O$
that agrees with $\gamma_n$ on $\partial\mathcal B$, the maps
$G\circ\widetilde\gamma$ and $G\circ\gamma_n$ coincide on the
boundary. Their Brouwer degrees are equal (use the straight homotopy
in the target $\R^4$, whose boundary values remain fixed).
Hence some point of $\widetilde\gamma(\mathcal B)$ belongs to
$\Neh_{\lambda,\nu}$, and
\begin{equation}\label{PS-topological-barrier}
 \max_{\xi\in\mathcal B}\A(\widetilde\gamma(\xi))\geq m.
\end{equation}

\emph{3. Local Lipschitz regularity needed for an actual flow.}
The functional $\A$ is $C^1$ on $\X$ and its first derivative is
Lipschitz on bounded subsets of $\X$, uniformly for $\lambda\in K$
and small $\nu$. To verify the latter assertion, for $2<s\leq2^*$
use the pointwise estimate
\[
 \big||a|^{s-2}a-|b|^{s-2}b\big|
 \leq C_s(|a|+|b|)^{s-2}|a-b|
\]
and H\"older's inequality in $L^s$. For the mixed integrand
$H(a,b)=|a|^\alpha|b|^\beta$, the assumptions
$\alpha,\beta\geq2$ imply
\[
 |\nabla H(\zeta)-\nabla H(\omega)|
 \leq C_{\alpha,\beta}(|\zeta|+|\omega|)^{r-2}|\zeta-\omega|,
 \qquad\zeta,\omega\in\R^2.
\]
Applying H\"older's inequality in $L^r$, and then the Sobolev
embeddings, yields for each $R>0$ a constant $L_R$ such that
\begin{equation}\label{PS-Lipschitz}
 \|\A'(x)-\A'(y)\|_{\X^*}\leq L_R\|x-y\|_\X
 \qquad(\|x\|_\X,\|y\|_\X\leq R).
\end{equation}
Here the Hilbert-space Riesz gradient $\nabla_\X\A$ is defined by
$\langle\A'(y),h\rangle=(\nabla_\X\A(y),h)_\X$ for all $h\in\X$.
Thus $\frac{\nabla_\X\A}{\|\nabla_\X\A\|_\X}$ is locally Lipschitz
where its denominator is positive.

\emph{4. Quantitative localized deformation and its contradiction.}
Fix $n$ so large that $\varepsilon:=\varepsilon_n$ satisfies
$\varepsilon<d_0$ and~\eqref{PS-boundary}.
Suppose, contrary to the conclusion, that there is no
$y\in\X$ satisfying the three inequalities
\begin{equation}\label{PS-negated}
 |\A(y)-m|\leq2\varepsilon^2,\qquad
 \dist(y,\Gamma_n)\leq2\varepsilon,\qquad
 \|\A'(y)\|_{\X^*}\leq8\varepsilon.
\end{equation}
Define the open energy--distance band
\[
 U=\{y\in\X:|\A(y)-m|<2\varepsilon^2,
                     \ \dist(y,\Gamma_n)<2\varepsilon\}.
\]
Then $\|\A'(y)\|_{\X^*}>8\varepsilon$ throughout $U$.
Choose Lipschitz cutoffs $\chi:\R\to[0,1]$ and
$\theta:[0,\infty)\to[0,1]$ such that
\[
 \begin{array}{ll}
 \chi(t)=1\ (|t|\leq1),&\chi(t)=0\ (|t|\geq2),\\
 \theta(t)=1\ (0\leq t\leq1),&\theta(t)=0\ (t\geq2).
 \end{array}
\]
On $U$ define
\[
 W(y)=-\chi\!\left(\frac{\A(y)-m}{\varepsilon^2}\right)
         \theta\!\left(\frac{\dist(y,\Gamma_n)}{\varepsilon}\right)
         \frac{\nabla_\X\A(y)}{\|\nabla_\X\A(y)\|_\X},
\]
and put $W(y)=0$ outside $U$.
By~\eqref{PS-Lipschitz}, $W$ is locally Lipschitz;
its defining cutoffs vanish on the boundary of $U$, where the
extension by zero is continuous and locally Lipschitz. Moreover
$\|W(y)\|_\X\leq1$. The ordinary differential equation
\begin{equation}\label{PS-flow}
 \dot\Phi(t,y)=W(\Phi(t,y)),\qquad\Phi(0,y)=y
\end{equation}
therefore has a global flow for finite time, with
\begin{equation}\label{PS-speed}
 \|\Phi(t,y)-y\|_\X\leq t,\qquad
 \frac{d}{dt}\A(\Phi(t,y))
 =-\chi\theta\,\|\A'(\Phi(t,y))\|_{\X^*}\leq0.
\end{equation}
We obtain the derivative identity from the chain rule and the Riesz
representation; the arguments of $\chi,\theta$ are as above.

Take the fixed deformation time $T=\frac{\varepsilon}{2}$ and let
$\widetilde\gamma_n(\xi)=\Phi(T,\gamma_n(\xi))$.
Using the speed bound, we obtain
$\dist(\Phi(t,\gamma_n(\xi)),\Gamma_n)\leq t\leq\frac{\varepsilon}{2}$
for $0\leq t\leq T$. Hence $\theta=1$ throughout these trajectories.
For as long as the energy lies in
$[m-\varepsilon^2,m+\varepsilon^2]$, we also have $\chi=1$,
and thus its rate of decrease is strictly less than
$-8\varepsilon$. Since
$\A(\gamma_n(\xi))\leq m+\varepsilon^2$, a trajectory starting
above $m-\varepsilon^2$ must cross that level within time
$\frac{2\varepsilon^2}{8\varepsilon}=\frac{\varepsilon}{4}<T$.
A trajectory starting below that level remains below it, by
monotonicity. We have proved the explicit upper bound
\begin{equation}\label{PS-lowered}
 \max_{\xi\in\mathcal B}\A(\widetilde\gamma_n(\xi))
 \leq m-\varepsilon^2.
\end{equation}
By~\eqref{PS-boundary}, the energy cutoff vanishes on
$\gamma_n(\partial\mathcal B)$, so the flow fixes the boundary.
For $0\leq t\leq T$, the images of the homotopy
$\xi\mapsto\Phi(t,\gamma_n(\xi))$ belong to $\mathcal O$
by~\eqref{O-tube} and the speed bound.
The degree argument~\eqref{PS-topological-barrier} consequently
forces the maximum in~\eqref{PS-lowered} to be at least $m$,
a contradiction. Therefore a point $y_n$ satisfying~\eqref{PS-negated}
exists for all sufficiently large $n$.

\emph{5. Quantitative proximity to the original minimizing sequence.}
Since $\Gamma_n$ is compact, choose $\xi_n\in\mathcal B$ with
\begin{equation}\label{PS-xi-near}
 \|y_n-\gamma_n(\xi_n)\|_\X
 =\dist(y_n,\Gamma_n)\leq2\varepsilon_n.
\end{equation}
All these points belong to a fixed ball in $\X$. By the
fundamental theorem of calculus,~\eqref{PS-Lipschitz},
$\|\A'(y_n)\|\leq8\varepsilon_n$, and~\eqref{PS-xi-near},
\begin{align}\label{PS-energy-near}
 |\A(\gamma_n(\xi_n))-\A(y_n)|
 &\leq\|\A'(y_n)\|\,\|\gamma_n(\xi_n)-y_n\|_\X
   +\frac{L_R}{2}\|\gamma_n(\xi_n)-y_n\|_\X^2\\
 &\leq(16+2L_R)\varepsilon_n^2.\nonumber
\end{align}
Combining~\eqref{chosen-min-seq},~\eqref{PS-negated}, and
\eqref{PS-energy-near}, we obtain
\[
 0\leq J_{z_n}(\one)-J_{z_n}(\xi_n)
 \leq (19+2L_R)\varepsilon_n^2.
\]
From the quadratic gap~\eqref{uniformquadraticgap}, we obtain
\[
 |\xi_n-\one|\leq\sqrt{\frac{19+2L_R}{c_{\mathrm{gap}}}}\,\varepsilon_n.
\]
The four sign parts of $z_n$ have uniformly bounded $H^1$ norms,
so $\|\gamma_n(\xi_n)-z_n\|_\X\leq C|\xi_n-\one|$.
Together with~\eqref{PS-xi-near}, this proves the last bound
in~\eqref{PSsequence}. The other bounds are exactly
\eqref{PS-negated}. The four sign parts of $y_n$ retain uniform
$L^p$ and $L^q$ lower bounds by truncation continuity and
$\|y_n-z_n\|_\X\to0$.

Finally, suppose a minimizer $z$ is already attained.
Apply the construction with the constant minimizing sequence
$z_n=z$ and any $\varepsilon_n\downarrow0$. We obtain
$y_n\to z$ in $\X$ and $\A'(y_n)\to0$ in $\X^*$.
Continuity of $\A'$ implies $\A'(z)=0$.
We have proved the assertion for \emph{every} minimizer without
assuming differentiability of the sign truncation maps.
\end{proof}

\section{Compactness and attainment}\label{sec:compact}
\begin{lemma}\label{lem:radialembedding}
For $2<s<2^*$, the embedding
\begin{equation}\label{radialcompact}
\Hr\hookrightarrow L^s(\RN),\qquad 2<s<2^*,
\end{equation}
is compact.
\end{lemma}
\begin{proof}
For radial $w\in H^1$, the estimate
$|w(x)|\leq C|x|^{-\frac{N-1}{2}}\|w\|_{H^1}$ holds almost everywhere. For smooth functions, we integrate the derivative of $\rho^{N-1}w(\rho)^2$ from $\rho=|x|$ to infinity. We obtain the general case by density. Thus, for $R>0$ and $s>2$,
\[
\int_{|x|>R}|w|^s\dx
\leq C R^{-\frac{(N-1)(s-2)}{2}}\|w\|_{H^1}^{s-2}\|w\|_2^2.
\]
The right-hand side tends to zero uniformly on bounded subsets of $H^1$. Applying Rellich's theorem on $\{|x|<R\}$ and then letting $R\to\infty$, we obtain~\eqref{radialcompact}; see also~\cite[Chapter 1]{Willem}. Under theassumption of~\eqref{H}, this applies to both power terms.
\end{proof}

\begin{proposition}\label{prop:signlimit}
Under theassumption of~\eqref{H}, fix $K$ and the action bound $L$ used in Section~\ref{sec:bounds}. Choose $\bar\nu>0$ so that we can apply Lemma~\ref{lem:nonvanish} uniformly for $0\leq\nu\leq\bar\nu$. Suppose
\[
\lambda_n\in K,\quad 0\leq\nu_n\leq\bar\nu,\quad
z_n=(u_n,v_n)\in\Neh_{\lambda_n,\nu_n},\quad
\A_{\lambda_n,\nu_n}(z_n)\leq L.
\]
Every weakly convergent subsequence has a limit $(u,v)$ satisfying
\begin{equation}\label{fourlimit}
\|u^\pm\|_p\geq\delta,\qquad \|v^\pm\|_q\geq\delta
\end{equation}
for a uniform $\delta>0$. In particular neither a semitrivial limit nor a limit with a sign-definite component can occur. The same conclusion holds for any sequence $y_n$ with $\|y_n-z_n\|_{H^1\times H^1}\to0$.
\end{proposition}
\begin{proof}
By~\eqref{coercive}, the sequence $z_n$ is bounded in $\X$, and by Lemma~\ref{lem:nonvanish}, all four  integrals are bounded away from zero. If $z_n\rightharpoonup(u,v)$, we use the compact embeddings~\eqref{radialcompact} to obtain
\[
u_n\to u\text{ in }L^p(\RN),\qquad v_n\to v\text{ in }L^q(\RN).
\]
For either sign, we use the pointwise truncation inequality to obtain
\[
\|u_n^\pm-u^\pm\|_p\leq\|u_n-u\|_p\to0,
\qquad \|v_n^\pm-v^\pm\|_q\leq\|v_n-v\|_q\to0.
\]
Thus the lower bounds pass through strong convergence:
$\|u^\pm\|_p=\lim_n\|u_n^\pm\|_p\geq\delta$ and
$\|v^\pm\|_q=\lim_n\|v_n^\pm\|_q\geq\delta$.
In particular, a zero component or a component of one sign contradicts one of these four inequalities. We have proved~\eqref{fourlimit}. We obtain the conclusion for $y_n$ by Sobolev embedding and the same truncation estimate.
\end{proof}

By Proposition~\ref{prop:signlimit}, the Palais--Smale sequence from Lemma~\ref{lem:PSsequence} has a nodal weak limit, including when $r=2^*$. We next rule out a norm defect.

\begin{lemma}\label{lem:subcompact}
Let $r<2^*$ and let $\lambda_n\to\lambda\in(0,\infty)^2$. Every bounded Palais--Smale sequence for $\A_{\lambda_n,\nu}|_\X$ has a strongly convergent subsequence in $H^1(\RN)^2$.
\end{lemma}
\begin{proof}
Let $z_n=(u_n,v_n)$ be such a sequence. After extracting a weakly convergent subsequence, we use~\eqref{radialcompact} to obtain strong convergence in $L^p\times L^q$ and in $L^r\times L^r$. The two mixed derivatives are
$\alpha|u|^{\alpha-2}u|v|^\beta$ and
$\beta|u|^\alpha|v|^{\beta-2}v$. For example, by H\"older's inequality,
\[
\bigl\||u|^{\alpha-2}u|v|^\beta\bigr\|_{\frac{r}{r-1}}
\leq\|u\|_r^{\alpha-1}\|v\|_r^\beta,
\]
because $(\alpha-1)+\beta=r-1$. By continuity of these power maps under strong $L^r$ convergence, both derivatives converge strongly in $L^{\frac{r}{r-1}}$. Pass to the limit against radial test functions to obtain the weak equation for the weak limit $(u,v)$. Subtract that equation from the equations for the sequence, and test with $(u_n-u,v_n-v)$. The power terms tend to zero by strong $L^p,L^q$ convergence, and the mixed terms tend to zero by the preceding estimate and H\"older's inequality. The terms containing $\lambda_n-\lambda$ tend to zero by the bounded $L^2$ norms. Since $\lambda_{n,j}$ are bounded away from zero, the full $H^1\times H^1$ norm of the difference tends to zero.
\end{proof}

For every $N\geq3$, define the sharp scalar Sobolev constant by
\begin{equation}\label{Sobconstant}
S=\inf_{0\ne\varphi\in C_c^\infty(\RN)}
\frac{\|\nabla\varphi\|_2^2}{\|\varphi\|_{2^*}^2}>0,
\end{equation}
where $C_c^\infty(\RN)$ denotes smooth compactly supported functions. The resulting inequality $S\|w\|_{2^*}^2\leq\|\nabla w\|_2^2$ holds for $w\in H^1$ by density. For the remainder of this section we assume $r=2^*$. For the interaction in~\eqref{mixedintegral}, put
\begin{equation}\label{Sab}
S_{\alpha\beta}=\frac{r}{\alpha^{\frac{\alpha}{r}}\beta^{\frac{\beta}{r}}}S.
\end{equation}
Combining the H\"older and Sobolev inequalities with the weighted arithmetic-geometric mean inequality, we obtain
\begin{equation}\label{vectorSob}
S_{\alpha\beta}C(u,v)^{\frac{2}{r}}
\leq\|\nabla u\|_2^2+\|\nabla v\|_2^2.
\end{equation}
Indeed, if $a_1=\|\nabla u\|_2^2$ and $a_2=\|\nabla v\|_2^2$, then
\[
C(u,v)^{\frac{2}{r}}\leq S^{-1}a_1^{\frac{\alpha}{r}}a_2^{\frac{\beta}{r}}
\leq S^{-1}(\frac{\alpha}{r})^{\frac{\alpha}{r}}(\frac{\beta}{r})^{\frac{\beta}{r}}(a_1+a_2).
\]
The critical concentration level is
\begin{align}
c_*(\nu)&=\frac1N S_{\alpha\beta}^{\frac{N}{2}}(\nu r)^{-\frac{N-2}{2}}\label{threshold}\\
&=\frac{2}{N-2}\nu^{-\frac{N-2}{2}}
\alpha^{-\frac{(N-2)\alpha}{4}}\beta^{-\frac{(N-2)\beta}{4}}S^{\frac{N}{2}}.
\nonumber
\end{align}
This is the mixed critical threshold appearing in~\cite{GHShZ}; in particular $c_*(\nu)\to\infty$ as $\nu\downarrow0$.

\begin{lemma}\label{lem:BL}
If $z_n=(u_n,v_n)$ is bounded in $L^r(\RN)^2$ and converges almost everywhere to $z=(u,v)$, then
\begin{equation}\label{BL}
C(u_n,v_n)=C(u,v)+C(u_n-u,v_n-v)+o(1).
\end{equation}
\end{lemma}
\begin{proof}
Define $F(s,t)=|s|^\alpha|t|^\beta$ for $(s,t)\in\R^2$. By homogeneity and $C^1$ regularity, we have $|\nabla F(\zeta)|\leq C|\zeta|^{r-1}$ for $\zeta\in\R^2$. Applying the fundamental theorem of calculus along a segment and then Young's inequality, we find that for every $\varepsilon>0$ there is $C_\varepsilon>0$ such that, for $\zeta,\omega\in\R^2$, 
\[
|F(\zeta+\omega)-F(\zeta)|\leq\varepsilon|\zeta|^r+C_\varepsilon|\omega|^r.
\]
Apply this to $\zeta=z_n-z$, $\omega=z$. The positive part of
$|F(z_n)-F(z_n-z)-F(z)|-\varepsilon|z_n-z|^r$
is bounded by an integrable multiple of $|z|^r$ and converges pointwise to zero. By dominated convergence, followed by $\varepsilon\downarrow0$, we obtain~\eqref{BL}.
\end{proof}

\begin{lemma}\label{lem:critcompact}
Suppose $r=2^*$, $\nu>0$, $\lambda_n\to\lambda\in(0,\infty)^2$, and $z_n=(u_n,v_n)\in\X$ is bounded with
\[
\A_{\lambda_n,\nu}(z_n)\to\ell<c_*(\nu),\qquad
(\A_{\lambda_n,\nu}|_\X)'(z_n)\to0.
\]
Then a subsequence converges strongly in $H^1(\RN)^2$.
\end{lemma}
\begin{proof}
Let $z_n\rightharpoonup z=(u,v)$ in $\X$. The pure nonlinearities converge strongly in their subcritical spaces by~\eqref{radialcompact}. The mixed derivatives converge almost everywhere and are bounded in $L^{\frac{r}{r-1}}$. Since $\frac{r}{r-1}>1$, we obtain uniform integrability on every bounded set and convergence of their integrals against smooth compactly supported test functions. Approximating a general $L^r$ test function by such tests and using the same uniform dual bound, we obtain weak convergence in $L^{\frac{r}{r-1}}$. Passing to the limit, we find that $z$ is a critical point of the restricted action.

Write $z_n-z=(\phi_n,\psi_n)$ and
\[
R_n=\|\phi_n\|_{1,\lambda_1}^2+\|\psi_n\|_{1,\lambda_2}^2.
\]
Testing the derivative condition at $z_n$ against $z_n$, we obtain the total Nehari identity with an $o(1)$ error. We subtract the corresponding identity at $z$. The quadratic terms split by weak convergence, the power integrals converge strongly, and we split the mixed term using~\eqref{BL}. We thus obtain the first identity below. For the second, we split the action and substitute the first identity, using $\frac{1}{2}-\frac{1}{r}=\frac{1}{N}$:
\begin{equation}\label{defect}
R_n=\nu rC(\phi_n,\psi_n)+o(1),\qquad
\ell=\A_{\lambda,\nu}(z)+\frac1N R_n+o(1).
\end{equation}
The terms involving $\lambda_n-\lambda$ tend to zero by boundedness. Pass to a subsequence for which $R_n$ has a limit. If this limit is positive, we use~\eqref{vectorSob} and the definition of $R_n$ to obtain
\[
R_n\geq S_{\alpha\beta}C(\phi_n,\psi_n)^{\frac{2}{r}}
=S_{\alpha\beta}\left(\frac{R_n}{\nu r}+o(1)\right)^{\frac{2}{r}}.
\]
Since $1-\frac{2}{r}=\frac{2}{N}$, we take the limit and raise to the power $\frac{N}{2}$ to obtain
$\lim \frac{R_n}{N}\geq S_{\alpha\beta}^{\frac{N}{2}}(\nu r)^{-\frac{N-2}{2}}\frac{1}{N}=c_*(\nu)$.
Every critical point, including a semitrivial one, has nonnegative action: testing against itself, we obtain
\[
\A_{\lambda,\nu}(z)=k_p\mu_1\|u\|_p^p+k_q\mu_2\|v\|_q^q
+(\frac{r}{2}-1)\nu C(u,v)\geq0.
\]
Combining this inequality with~\eqref{defect}, we contradict $\ell<c_*(\nu)$. Thus $R_n\to0$. Since both limiting frequencies are positive, we obtain strong convergence in the full $H^1\times H^1$ norm, including the $L^2$ parts.
\end{proof}

\begin{proof}[Attainment and the least-action identity in Theorem~\ref{thm:fixed}]
Use Lemma~\ref{lem:PSsequence}. In the subcritical case, apply Lemma~\ref{lem:subcompact}. In the critical case decrease $\nu_K$ so that $L<c_*(\nu)$ for $0<\nu<\nu_K$, and apply Lemma~\ref{lem:critcompact}. At $\nu=0$, the two equations decouple and only subcritical pure nonlinearities remain, so compactness follows directly from~\eqref{radialcompact}.

The strong limit has four nonzero sign parts and satisfies the four constraints. It attains $m_\nu(\lambda)$ and is a restricted critical point. Each attained minimizer is a restricted critical point by the last assertion of Lemma~\ref{lem:PSsequence}. To obtain the full weak equation, average an arbitrary test pair over the orthogonal group. At a radial point the action derivative on a test pair equals its value on the radial average, which is zero. Thus every minimizer is a full weak solution.

Conversely, every radial nodal solution belongs to~\eqref{Nnod}, by testing each equation against its positive and negative parts. We have proved~\eqref{least-solutions}. We obtain positivity of $m_\nu$ from~\eqref{coercive} and Lemma~\ref{lem:nonvanish}. We prove compactness along varying frequencies and continuity of the level next.
\end{proof}

\begin{corollary}\label{cor:constructionlimit}
Assume~\eqref{H}, and keep the compact positive frequency rectangle $K$ and the action bound $L$ from Section~\ref{sec:bounds}. Choose $\bar\nu>0$ so that we can apply Lemma~\ref{lem:nonvanish} for $0\leq\nu\leq\bar\nu$; when $r=2^*$, also require $L<c_*(\bar\nu)$. Suppose
\[
\lambda_n\in K,\quad\lambda_n\to\lambda,\qquad
0\leq\nu_n\leq\bar\nu,\quad\nu_n\to\nu,
\]
and let $z_n\in\Neh_{\lambda_n,\nu_n}$ satisfy
$\A_{\lambda_n,\nu_n}(z_n)\leq L$. If $y_n\in\X$ satisfies
\[
\|y_n-z_n\|_{H^1\times H^1}\to0,
\qquad (\A_{\lambda_n,\nu_n}|_\X)'(y_n)\to0,
\]
then, along a subsequence, both $y_n$ and $z_n$ converge strongly in $H^1(\RN)^2$ to a full weak solution $(u,v)$ at $(\lambda,\nu)$. All four sign parts of this limit are nonzero. If $y_n=(\widetilde u_n,\widetilde v_n)$ has masses tending to $a,b$, respectively, then $\|u\|_2^2=a$ and $\|v\|_2^2=b$.
\end{corollary}
\begin{proof}
By coercivity, $z_n$, and hence $y_n$, are bounded in $\X$. Using Proposition~\ref{prop:signlimit}, we exclude trivial, semitrivial and seminodal weak limits before considering the critical norm defect. On this bounded set, the mixed derivative has a uniformly bounded dual norm by H\"older's inequality. Consequently
\[
\bigl\|(\A_{\lambda,\nu}|_\X)'(y_n)
-(\A_{\lambda_n,\nu_n}|_\X)'(y_n)\bigr\|
\leq C\bigl(|\lambda_n-\lambda|+|\nu_n-\nu|\bigr)\to0.
\]
The actions differ by a quantity tending to zero as well. After extracting an action-convergent subsequence, its level is at most $L$.

If $r<2^*$, apply Lemma~\ref{lem:subcompact} at the limiting parameters. If $r=2^*$ and $\nu>0$, the level is less than $c_*(\bar\nu)\leq c_*(\nu)$, so we may apply Lemma~\ref{lem:critcompact}. If $r=2^*$ and $\nu=0$, the mixed derivative disappears and the limiting functional contains only the two subcritical powers. Using strong $L^p,L^q$ convergence and testing the difference of the equations, we obtain strong $H^1\times H^1$ convergence, because both limiting frequencies are positive. In all three cases, we recover the full weak equations by radial averaging.

By strong $H^1$ convergence, we also obtain strong $L^2$ convergence. Thus
$\|u\|_2^2=\lim_n\|\widetilde u_n\|_2^2$ and
$\|v\|_2^2=\lim_n\|\widetilde v_n\|_2^2$.
The four nonzero signs follow from~\eqref{fourlimit}.
\end{proof}

\section{The level as a function of two frequencies}\label{sec:parameter}
We begin with an estimate for the action level as the coupling tends to zero.

\begin{proposition}\label{prop:limit}
The scalar radial nodal minima $d_{p,\mu_1}$ and $d_{q,\mu_2}$ are attained and their minimizing sets are compact. Furthermore,
\begin{equation}\label{mzero}
m_0(\lambda)=d_{p,\mu_1}(\lambda_1)+d_{q,\mu_2}(\lambda_2),
\end{equation}
and there is a constant $C_K>0$, independent of $\lambda\in K$ and small $\nu$, such that
\begin{equation}\label{Oestimate}
|m_\nu(\lambda)-m_0(\lambda)|\leq C_K\nu.
\end{equation}
If $\nu_n\downarrow0$, $\lambda_n\to\lambda\in K$, and $z_n$ minimizes $m_{\nu_n}(\lambda_n)$, then a subsequence converges strongly to a product of scalar radial nodal action ground states at $\lambda$.
\end{proposition}
\begin{proof}
We first prove the scalar assertion for every $N\geq3$, including the dimensions excluded by~\eqref{H}. Fix $2<s<2^*$, $\mu>0$ and $\ell>0$. A radial function with two nonzero sign parts is projected to the scalar nodal constraints by the unique amplitudes
\[
t_\pm=\left(\frac{\|w^\pm\|_{1,\ell}^2}
{\mu\|w^\pm\|_s^s}\right)^{\frac{1}{s-2}}.
\]
We have proved nonemptiness. On those constraints,
$I_{s,\mu,\ell}(w)=k_s\|w\|_{1,\ell}^2$, and
\[
\|w^\pm\|_{1,\ell}^2=\mu\|w^\pm\|_s^s
\leq C\bigl(\|w^\pm\|_{1,\ell}^2\bigr)^{\frac{s}{2}}.
\]
Hence bounded-action sequences are bounded in $H^1$, with both sign norms bounded away from zero. Along $(t_+,t_-)\mapsto t_+w^+-t_-w^-$ the action separates into two scalar functions. Its unique maximum is at $(1,1)$, with diagonal Hessian entries $-(s-2)\mu\|w^\pm\|_s^s<0$. We may therefore use the deformation argument from Lemma~\ref{lem:PSsequence} with two amplitudes: the determinant of this two-dimensional Hessian is positive, the degree is nonzero, and a bounded Palais--Smale sequence lies at distance tending to zero from a nodal minimizing sequence. This uses no interaction exponent or dimension restriction from~\eqref{H}.

By radial compactness, we obtain strong $L^s$ convergence. Passing to the scalar weak equation and testing the difference, we also obtain strong $H^1$ convergence, since $\ell>0$. Both signs survive by the $L^s$ lower bounds. Thus we attain the infimum and use the same deformation argument to show that every minimizer is a radial critical point. By averaging test functions over rotations, we see that it is a full scalar weak solution. The estimates are uniform on compact positive frequency intervals; applying the same argument to minimizers, we obtain compactness along convergent frequencies.

When $\nu=0$, the vector admissible set is the product of the two scalar nodal sets. Taking the infimum, we obtain~\eqref{mzero}.

Starting from a product of scalar minimizers, we use the finite-dimensional implicit-function theorem to project the product onto~\eqref{Nnod} by adjusting its four sign amplitudes. The coefficients and the inverse diagonal derivative are uniformly bounded on $K$. The amplitudes equal $1+O(\nu)$, so the change of action is $O(\nu)$. We obtain $m_\nu\leq m_0+C_K\nu$.

For the reverse inequality take an attained minimizer $z=(u,v)$ at coupling $\nu$. Its four sign parts satisfy~\eqref{uniformcoeff}. Set
\begin{equation}\label{projectzero}
s_\sigma=\left(\frac{A_\sigma}{B_\sigma}\right)^{\frac{1}{p-2}},
\qquad t_\tau=\left(\frac{D_\tau}{E_\tau}\right)^{\frac{1}{q-2}}.
\end{equation}
Then $\widehat z=\gamma_z(s_+,s_-,t_+,t_-)$ lies on the decoupled nodal Nehari set. By~\eqref{coeffconstraints} and~\eqref{uniformcoeff}, all these amplitudes are $1+O(\nu)$. Using uniform boundedness, we find
\[
m_0(\lambda)\leq\A_{\lambda,0}(\widehat z)
\leq\A_{\lambda,\nu}(z)+C_K\nu.
\]
We have proved~\eqref{Oestimate}.

For the last statement, $z_n$ is bounded and solves the system. On bounded subsets of $H^1\times H^1$, the derivative of the mixed integral is bounded in the dual, including when $r=2^*$. Since $\nu_n\to0$, the sequence is Palais--Smale for the decoupled subcritical action at the limiting frequencies. By radial compactness, we obtain strong convergence. The four sign parts remain nonzero, and from~\eqref{Oestimate} we identify the limit as a minimizer of~\eqref{mzero}.
\end{proof}

\begin{proposition}\label{prop:envelope}
Fix $0\leq\nu<\nu_K$. The level $m_\nu$ is continuous on $K$ and locally Lipschitz on $\operatorname{int} K$. Minimizers are compact along frequency sequences in $K$. If $z=(u,v)$ minimizes $m_\nu(\lambda)$ and $h=(h_1,h_2)\in\R^2$ is a frequency increment, then, as $h\to0$ with $\lambda+h\in K$,
\begin{equation}\label{envelope}
m_\nu(\lambda+h)\leq m_\nu(\lambda)
+\frac12h_1\|u\|_2^2+\frac12h_2\|v\|_2^2+O(|h|^2).
\end{equation}
The constants can be chosen uniformly on smaller positive frequency rectangles.
\end{proposition}
\begin{proof}
In Section~\ref{sec:compact}, we obtained a minimizer at every frequency in $K$ and proved that every such minimizer is a restricted critical point. Its action is at most $L_0$, so the coefficient bounds~\eqref{uniformcoeff} and the Hessian estimate~\eqref{Hnegative} apply uniformly.

Fix one minimizer $z=(u,v)$ at $\lambda$. Let $G_{\lambda+h}$ denote the vector of the four constraint functions at the perturbed frequency. By~\eqref{G-on-path},
\[
 D_\xi\bigl[G_{\lambda+h}(\gamma_z(\xi))\bigr]_{(h,\xi)=(0,\one)}
 =D^2J_z(\one).
\]
The inverse matrix has norm at most $\frac{1}{\eta}$. The coefficients in the four constraints, their derivatives in positive amplitudes near $\one$, and the sign masses are uniformly bounded. We therefore apply the finite-dimensional implicit-function theorem to obtain positive amplitudes $\xi(h)$ with
\[
 G_{\lambda+h}(\gamma_z(\xi(h)))=0,
 \qquad |\xi(h)-\one|\leq C|h|,
\]
on a neighborhood independent of the chosen minimizer and of the parameters in the stated compact ranges. Thus $\gamma_z(\xi(h))$ is admissible at $\lambda+h$ and retains all four sign parts.

The derivative of $J_z$ at $\one$ vanishes by the four constraints. By Taylor expansion in the amplitudes, we obtain
\[
 \A_{\lambda,\nu}(\gamma_z(\xi(h)))
 =\A_{\lambda,\nu}(z)+O(|h|^2).
\]
Each component mass changes by $O(|h|)$ under this amplitude adjustment. Since
\[
 \A_{\lambda+h,\nu}(w_1,w_2)-\A_{\lambda,\nu}(w_1,w_2)
 =\frac12h_1\|w_1\|_2^2+\frac12h_2\|w_2\|_2^2,
\]
we obtain
\[
 \A_{\lambda+h,\nu}(\gamma_z(\xi(h)))
 =m_\nu(\lambda)+\frac12h_1\|u\|_2^2
  +\frac12h_2\|v\|_2^2+O(|h|^2).
\]
Taking the infimum at $\lambda+h$, we prove~\eqref{envelope} for each minimizer separately and obtain upper semicontinuity. These are finite-dimensional derivatives of positive amplitudes.

For the reverse semicontinuity, let $\lambda_n\to\lambda$ in $K$ and let $z_n$ minimize $m_\nu(\lambda_n)$. From the preceding attainment argument, we obtain
\[
 (\A_{\lambda_n,\nu}|_\X)'(z_n)=0,
 \qquad \A_{\lambda_n,\nu}(z_n)\leq L_0<L.
\]
The sequence is bounded in $\X$. At the limiting frequency,
\[
 \|(\A_{\lambda,\nu}|_\X)'(z_n)\|_{\X^*}
 \leq C|\lambda_n-\lambda|,
 \qquad
 |\A_{\lambda,\nu}(z_n)-m_\nu(\lambda_n)|
 \leq C|\lambda_n-\lambda|.
\]
We may apply Corollary~\ref{cor:constructionlimit} with $y_n=z_n$ and constant coupling. In the critical case its uniform condition $L<c_*(\bar\nu)$ is retained in the choice of the common coupling bound. Along a subsequence, $z_n\to z$ strongly in $\X$; the limit retains all four signs and belongs to $\Neh_{\lambda,\nu}$. Hence
\[
 m_\nu(\lambda)\leq\A_{\lambda,\nu}(z)
 =\lim_n m_\nu(\lambda_n)\leq m_\nu(\lambda),
\]
where the last inequality follows from upper semicontinuity. Applying the same reasoning to every subsequence, we obtain continuity and see that every limit of minimizers is a minimizer at $\lambda$. We have proved compactness along frequency sequences.

Finally, apply~\eqref{envelope} at minimizers for each of two neighboring frequencies. Using the uniform mass bounds and amplitude estimates, we compare both directions and obtain a local Lipschitz bound.
\end{proof}

For the scalar laws used below, fix $2<s<2^*$ and $\mu>0$. For $\ell>0$, define the scalar dilation
\begin{equation}\label{scalardilation}
w_\ell(x)=\ell^{\frac{1}{s-2}}w(\sqrt\ell\,x).
\end{equation}
With the change of variables $y=\sqrt\ell\,x$, we obtain
\begin{align*}
\|w_\ell\|_2^2&=\ell^{\frac{2}{s-2}-\frac{N}{2}}\|w\|_2^2,\\
\|\nabla w_\ell\|_2^2&=\ell^{1+\frac{2}{s-2}-\frac{N}{2}}\|\nabla w\|_2^2,\\
\|w_\ell\|_s^s&=\ell^{\frac{s}{s-2}-\frac{N}{2}}\|w\|_s^s
=\ell^{1+\frac{2}{s-2}-\frac{N}{2}}\|w\|_s^s.
\end{align*}
The same identities hold for each sign part. Thus we map the scalar radial nodal Nehari set at frequency $1$ bijectively onto that at $\ell$ by dilation and multiply the action by $\ell^{1+\frac{2}{s-2}-\frac{N}{2}}$. With~\eqref{sigmakappa},
\begin{equation}\label{scalarlaw}
m_0(\lambda)=d_1\lambda_1^{\sigma_1+1}
+d_2\lambda_2^{\sigma_2+1}.
\end{equation}
Every scalar radial nodal action minimizer $w$ for exponent $s$, coefficient $\mu$ and frequency $\ell$ satisfies
\begin{equation}\label{masslaw}
\|w\|_2^2=2d_{s,\mu}'(\ell)
=2\left(1+\frac{2}{s-2}-\frac N2\right)
 d_{s,\mu}(1)\ell^{\frac{2}{s-2}-\frac{N}{2}}.
\end{equation}
Indeed, for a scalar minimizer $w$ at $\ell$, the amplitudes
\[
 t_\pm(h)=\left(1+
 \frac{h\|w^\pm\|_2^2}{\|w^\pm\|_{1,\ell}^2}\right)^{\frac{1}{s-2}}
 =1+O(|h|)
\]
project its sign parts onto the scalar nodal constraints at $\ell+h$. Expanding the action at these amplitudes gives
\[
 d_{s,\mu}(\ell+h)\leq d_{s,\mu}(\ell)
 +\frac h2\|w\|_2^2+O(h^2).
\]
From the dilation formula, we see that $d_{s,\mu}$ is differentiable. Comparing positive and negative increments in the last inequality, we obtain $\|w\|_2^2=2d_{s,\mu}'(\ell)$ for every minimizer and prove~\eqref{masslaw}.

\begin{lemma}\label{lem:mass}
Suppose $\lambda\in\operatorname{int} K$ is a local minimum of $F_{\nu,a,b}$ from~\eqref{Fdef}. Every minimizer $z=(u,v)$ of $m_\nu(\lambda)$ satisfies
$\|u\|_2^2=a$ and $\|v\|_2^2=b$. Moreover~\eqref{gradient} holds at $\lambda$.
\end{lemma}
\begin{proof}
From local minimality and~\eqref{envelope}, we obtain, for small $h$,
\begin{equation}\label{squeeze}
0\leq F_{\nu,a,b}(\lambda+h)-F_{\nu,a,b}(\lambda)
\leq\frac12h_1(\|u\|_2^2-a)
+\frac12h_2(\|v\|_2^2-b)+O(|h|^2).
\end{equation}
Let $e_1=(1,0)$ and $e_2=(0,1)$ be the two coordinate vectors. For $t>0$, substitute $h=te_1$ and $h=-te_1$ into~\eqref{squeeze}. After division by $t$,
\[
0\leq\frac12(\|u\|_2^2-a)+O(t),\qquad
0\leq-\frac12(\|u\|_2^2-a)+O(t).
\]
Letting $t\downarrow0$ in the two inequalities, we find $\|u\|_2^2=a$. Taking $h=\pm te_2$ in the same way, we obtain $\|v\|_2^2=b$. The minimizer $z$ was arbitrary. Substituting these equalities back into~\eqref{squeeze}, we obtain
\[
0\leq m_\nu(\lambda+h)-m_\nu(\lambda)-\frac12(a,b)\cdot h
\leq C|h|^2,
\]
which is the claimed differentiability and gradient identity.
\end{proof}

For the mass selection in Lemma~\ref{lem:mass}, we use a two-frequency version of the one-sided comparison in~\cite[Lemma 4.2 and Remark 4.3]{DCDGS}.

\section{Proof of the normalization theorem}\label{sec:normalization}
\begin{proof}[Proof of Theorem~\ref{thm:mass}]
For the two exponents $p,q$, formula~\eqref{masslaw} becomes $2(\sigma_j+1)d_j\ell^{\sigma_j}$, $j=1,2$. Under~\eqref{M}, both expressions increase continuously from zero to infinity. Choose positive endpoints $\ell_j<L_j$, $j=1,2$, and the rectangle
\[
K=[\ell_1,L_1]\times[\ell_2,L_2]\subset(0,\infty)^2
\]
so that
\begin{equation}\label{Kchoice}
\begin{aligned}
2(\sigma_1+1)d_1\ell_1^{\sigma_1}&<a_-,&
2(\sigma_1+1)d_1L_1^{\sigma_1}&>a_+,\\
2(\sigma_2+1)d_2\ell_2^{\sigma_2}&<b_-,&
2(\sigma_2+1)d_2L_2^{\sigma_2}&>b_+.
\end{aligned}
\end{equation}
After fixing $K$, we choose the coupling bound so that the uniform estimates and compactness conclusions of Sections~\ref{sec:bounds}--\ref{sec:parameter} apply. When $r=2^*$, we retain $L<c_*(\bar\nu)$ from Corollary~\ref{cor:constructionlimit} and require $\nu<\bar\nu$. All further restrictions below are uniform for $(a,b)\in\mathcal M$.

For each target mass pair, we find from~\eqref{scalarlaw} that
\[
F_{0,a,b}(\lambda)=d_1\lambda_1^{\sigma_1+1}
+d_2\lambda_2^{\sigma_2+1}-\frac12(a\lambda_1+b\lambda_2).
\]
It is strictly convex and has the unique minimizer $\lambda^0$ in~\eqref{lambda-limit}. By~\eqref{Kchoice} this minimizer belongs to $\operatorname{int} K$. Its distance from $\partial K$ is uniformly positive for $(a,b)\in\mathcal M$. The Hessian is diagonal, with entries
$d_j\sigma_j(\sigma_j+1)\lambda_j^{\sigma_j-1}$, $j=1,2$. They have a positive minimum on $K$, so the Hessian is bounded below by $cI_2$, where $I_2$ is the $2\times2$ identity matrix and $c>0$ is independent of the target masses. Consequently there is a uniform boundary gap
\begin{equation}\label{frequencygap}
\min_{\lambda\in\partial K}F_{0,a,b}(\lambda)
\geq F_{0,a,b}(\lambda^0)+\delta,\qquad\delta>0.
\end{equation}

By Theorem~\ref{thm:fixed} and Proposition~\ref{prop:envelope}, $F_{\nu,a,b}$ is continuous on $K$ and therefore attains its minimum. For $\lambda\in\partial K$, we use~\eqref{Oestimate} to obtain
\[
F_{\nu,a,b}(\lambda)
\geq F_{0,a,b}(\lambda^0)+\delta-C_K\nu
\geq F_{\nu,a,b}(\lambda^0)+\delta-2C_K\nu.
\]
If $C_K\nu<\frac{\delta}{3}$, the right-hand side is strictly greater than $F_{\nu,a,b}(\lambda^0)$. Every minimum on $K$ therefore lies in $\operatorname{int} K$.

We next estimate the distance between a minimizing frequency and $\lambda^0$. If $\lambda_\nu$ is any minimizer, the uniform Hessian lower bound and $F_{\nu,a,b}(\lambda_\nu)\leq F_{\nu,a,b}(\lambda^0)$ give
\[
 \frac c2|\lambda_\nu-\lambda^0|^2
 \leq F_{0,a,b}(\lambda_\nu)-F_{0,a,b}(\lambda^0)
 \leq2C_K\nu.
\]
Thus $|\lambda_\nu-\lambda^0|\leq C\sqrt\nu$ uniformly on $\mathcal M$, which proves~\eqref{lambda-limit}. By Lemma~\ref{lem:mass}, every action minimizer at such a frequency has the prescribed masses. We obtain the equations and the least-action property from Theorem~\ref{thm:fixed}.

For fixed target masses, we apply Corollary~\ref{cor:constructionlimit} with $y_n=z_n=(u_{\nu_n},v_{\nu_n})$, where $\nu_n\downarrow0$ and the frequencies tend to $\lambda^0\in(0,\infty)^2$. The action derivatives vanish because these are attained minimizers. We obtain strong $H^1\times H^1$ convergence along subsequences, and the limit satisfies
\[
\|u_0\|_2^2=a>0,\qquad \|v_0\|_2^2=b>0,
\qquad u_0^\pm\ne0,\quad v_0^\pm\ne0.
\]
Thus the limit $(u_0,v_0)$ is neither trivial nor semitrivial. By Proposition~\ref{prop:limit}, its components are scalar radial nodal action ground states. Since each original state has masses $(a,b)$,
\begin{equation}\label{EminF}
\E_\nu(u_\nu,v_\nu)=F_{\nu,a,b}(\lambda_\nu)=\min_K F_{\nu,a,b}.
\end{equation}
Using $a=2(\sigma_1+1)d_1(\lambda_1^0)^{\sigma_1}$ and the analogous formula for $b$, the limit of~\eqref{EminF} is exactly~\eqref{negative-limit}.
\end{proof}

\section{Critical powers and subcritical coupling}\label{sec:criticalpure}
Here $N\geq3$, the two powers are equal to an exponent $s$, and the assumptions are
\begin{equation}\label{reverse-assumptions}
p=q=s,\qquad \alpha,\beta>1,\qquad
2<r=\alpha+\beta<p_c=2+\frac4N,
\qquad s\in\{p_c,2^*\}.
\end{equation}
We retain radial symmetry and use~\eqref{reverse-assumptions} in place of~\eqref{H}--\eqref{M}. Under a coefficient-balance condition introduced below, we reduce the system to a scalar equation and prescribe masses on a ray. We minimize the action within the corresponding proportionality class. For comparison, positive normalized solutions with critical powers and subcritical coupling are known from~\cite{ZZZ}. Define the exponent of the mixed integral under mass-preserving dilation by
\begin{equation}\label{deltar}
\delta_r=\frac{N(r-2)}2.
\end{equation}
Thus $\delta_r<2$ under~\eqref{reverse-assumptions}.

\subsection{Exact reduction in a proportionality class}
For either value of $s$ in~\eqref{reverse-assumptions}, assume
\begin{equation}\label{balance}
\mu_1\alpha^{\frac{s-2}{2}}=\mu_2\beta^{\frac{s-2}{2}}=:\bar\mu>0,
\qquad \kappa:=\nu\alpha^{\frac{\alpha}{2}}\beta^{\frac{\beta}{2}}>0.
\end{equation}
Here $\bar\mu$ and $\kappa$ are the reduced power and coupling coefficients. Define
\begin{equation}\label{proportionalspace}
\mathcal Y_{\alpha\beta}
=\{(\sqrt\alpha\,w,\sqrt\beta\,w):w\in\Hr\},
\qquad \lambda_1=\lambda_2=\lambda>0.
\end{equation}
For $w\in\Hr$, define the reduced scalar action
\begin{equation}\label{reducedaction}
I_{\lambda,\bar\mu,\kappa}(w)
=\frac12\|w\|_{1,\lambda}^2
-\frac{\bar\mu}{s}\|w\|_s^s-\frac\kappa r\|w\|_r^r.
\end{equation}
\begin{lemma}\label{lem:proportional-reduction}
Assume~\eqref{reverse-assumptions} and~\eqref{balance}. A pair $(\sqrt\alpha\,w,\sqrt\beta\,w)\in\mathcal Y_{\alpha\beta}$ solves the first two equations in Equations~\eqref{system} at frequencies $(\lambda,\lambda)$ if and only if $w$ solves
\begin{equation}\label{reducedscalar}
-\Delta w+\lambda w=\bar\mu|w|^{s-2}w+\kappa|w|^{r-2}w.
\end{equation}
The component masses are $(\alpha\|w\|_2^2,\beta\|w\|_2^2)$, and
\begin{equation}\label{actionreduction}
\A_{(\lambda,\lambda),\nu}(\sqrt\alpha\,w,\sqrt\beta\,w)
=rI_{\lambda,\bar\mu,\kappa}(w).
\end{equation}
The four nodal constraints reduce to
\begin{equation}\label{scalarmixedN}
w^\pm\ne0,\qquad
\|w^\pm\|_{1,\lambda}^2
=\bar\mu\|w^\pm\|_s^s+\kappa\|w^\pm\|_r^r.
\end{equation}
A scalar radial nodal action minimizer therefore gives a least-action nodal solution within $\mathcal Y_{\alpha\beta}$.
\end{lemma}
\begin{proof}
Substitute $u=\sqrt\alpha\,w$ and $v=\sqrt\beta\,w$ into the first equation and divide by $\sqrt\alpha$. Its power coefficient becomes $\mu_1\alpha^{\frac{s-2}{2}}=\bar\mu$, and its mixed coefficient becomes
\[
\frac{\nu\alpha(\sqrt\alpha)^{\alpha-1}(\sqrt\beta)^\beta}
{\sqrt\alpha}
=\nu\alpha^{\frac{\alpha}{2}}\beta^{\frac{\beta}{2}}=\kappa.
\]
Dividing the second equation by $\sqrt\beta$, we find the same coefficients, which proves the equivalence with~\eqref{reducedscalar}. The mass identities follow directly from the substitution. Since $\mu_1\alpha^{\frac{s}{2}}+\mu_2\beta^{\frac{s}{2}}=r\bar\mu$ and $\alpha+\beta=r$, substitution into the vector action gives~\eqref{actionreduction}. For each sign $\sigma\in\{+,-\}$, the first-component constraint is $\alpha$ times
$\|w^\sigma\|_{1,\lambda}^2-\bar\mu\|w^\sigma\|_s^s-\kappa\|w^\sigma\|_r^r$, and the second-component constraint is $\beta$ times the same expression. Thus the four constraints reduce to~\eqref{scalarmixedN}. The asserted least-action property follows from~\eqref{actionreduction}.
\end{proof}

\subsection{Sobolev-critical powers: scalar construction and mass range}
\begin{lemma}\label{lem:scalarunit}
Let $N\geq3$, $s=2^*$, and $2<r<p_c$. There exists $\varepsilon_*(N,r)>0$ such that, for every $0\leq\varepsilon<\varepsilon_*$, some $\omega>0$ admits a radial nodal action ground state $W$ satisfying
\begin{equation}\label{unitproblem}
-\Delta W+\omega W=|W|^{r-2}W+\varepsilon|W|^{s-2}W,
\qquad \|W\|_2^2=1.
\end{equation}
Its energy is negative. The action comparison is over the entire scalar radial nodal Nehari set at $\omega$.
\end{lemma}
\begin{proof}
\emph{The unperturbed frequency and uniform action bound.}
Let $d_r=d_{r,1}(1)>0$ be the pure-$r$ scalar radial nodal level, whose attainment in every dimension $N\geq3$ was proved in Proposition~\ref{prop:limit}. Define
\[
\sigma=\frac{2}{r-2}-\frac N2>0,\quad
c_0(\omega)=d_r\omega^{\sigma+1},\quad
\omega_0=[2(\sigma+1)d_r]^{-\frac{1}{\sigma}},\quad
K_r=[\frac{\omega_0}{2},2\omega_0].
\]
Here $c_0$ is the unperturbed action level, and $\omega_0$ is the unique minimizer of $c_0(\omega)-\frac{\omega}{2}$ on $(0,\infty)$. For $\omega\in K_r$, define
\[
J_{\omega,\varepsilon}(W)
=\frac12\|W\|_{1,\omega}^2-\frac1r\|W\|_r^r
-\frac\varepsilon s\|W\|_s^s,
\]
and let $c_\varepsilon(\omega)$ be its infimum over radial $W$ with $W^\pm\ne0$ and
$\|W^\pm\|_{1,\omega}^2=\|W^\pm\|_r^r+\varepsilon\|W^\pm\|_s^s$.
Only in this proof, for a fixed nodal function $W$, abbreviate the quadratic and pure-$r$ sign terms by
\[
A_\pm=\|W^\pm\|_{1,\omega}^2,\qquad
B_\pm=\|W^\pm\|_r^r.
\]
Each amplitude $t>0$ is projected onto its constraint by solving
$A_\pm=B_\pm t^{r-2}+\varepsilon \|W^\pm\|_s^s t^{s-2}$. The right side is strictly increasing from zero to infinity, so the root is unique and is the strict maximum of the sign action. The two sign actions add. Testing with a pure-$r$ radial nodal minimizer and projecting its signs, we therefore obtain
\begin{equation}\label{scalarupper}
c_\varepsilon(\omega)\leq c_0(\omega),\qquad\omega\in K_r.
\end{equation}
Indeed, the mixed-power action is at most the pure-$r$ action at every amplitude, and the latter has its maximum at the unscaled pure-$r$ minimizer.

Recall $k_t=\frac{1}{2}-\frac{1}{t}$ for $t>2$. For this scalar proof, set $L=1+\max_{K_r}c_0$ and $M=\frac{L}{k_r}$. On the scalar nodal set of action at most $L$,
\[
J_{\omega,\varepsilon}(W)
=k_r\|W\|_r^r+k_s\varepsilon\|W\|_s^s
\geq k_r\|W\|_{1,\omega}^2,
\qquad \|W\|_{1,\omega}^2\leq M.
\]
For the inequality, we used $k_s>k_r$ and the sum of the two nodal constraints.

\emph{Preservation of both signs and exclusion of a critical residual.}
By Sobolev's inequality,
\[
\varepsilon \|W^\pm\|_s^s\leq\varepsilon S^{-\frac{s}{2}}A_\pm^{\frac{s}{2}}
\leq\varepsilon S^{-\frac{s}{2}}M^{\frac{s-2}{2}}A_\pm.
\]
Require
\begin{equation}\label{epsabs}
\varepsilon\leq\frac{S^{\frac{s}{2}}}{2M^{\frac{s-2}{2}}}.
\end{equation}
Then $B_\pm=A_\pm-\varepsilon \|W^\pm\|_s^s\geq \frac{A_\pm}{2}$. By the uniform embedding on $K_r$, we also have $B_\pm\leq C A_\pm^{\frac{r}{2}}$. Dividing by $A_\pm>0$, we obtain positive lower bounds for $A_\pm$ and $B_\pm$, uniformly in $\omega$ and the indicated $\varepsilon$ range.

At a nodal point, the two-amplitude Hessian is diagonal with entries
$-(r-2)B_\pm-(s-2)\varepsilon \|W^\pm\|_s^s$, which are uniformly negative. Using the unique sign maxima and the two-dimensional deformation argument from Proposition~\ref{prop:limit}, we construct a bounded Palais--Smale sequence $W_n$ at $c_\varepsilon(\omega)$, at distance tending to zero from a nodal minimizing sequence. Fix $\varepsilon$ and $\omega$ during this minimizing-sequence limit. Both $L^r$ sign norms remain bounded away from zero. Taking a weak limit $W$ and using radial $L^r$ compactness, we find
\begin{equation}\label{scalar-nonzerolimit}
\|(W_n)^\pm-W^\pm\|_r\to0,
\qquad \|W^\pm\|_r=\lim_n\|(W_n)^\pm\|_r\geq c>0.
\end{equation}
Thus the weak limit is nonzero and changes sign. Weak convergence of the critical derivative, justified by its almost everywhere convergence and bounded $L^{\frac{s}{s-1}}$ norm, gives $J_{\omega,\varepsilon}'(W)=0$ on the radial space.

Set $R_n=\|W_n-W\|_{1,\omega}^2$, the squared scalar residual norm. Applying the scalar version of the splitting argument in Lemma~\ref{lem:BL} to $|w|^s$ and using the derivative identities, we obtain
\begin{equation}\label{scalardefect}
R_n=\varepsilon\|W_n-W\|_s^s+o(1),\qquad
c_\varepsilon(\omega)=J_{\omega,\varepsilon}(W)+\frac{R_n}{N}+o(1).
\end{equation}
By testing at the weak critical point, we obtain $J_{\omega,\varepsilon}(W)\geq0$. If $R_n$ does not tend to zero, we use boundedness to extract a subsequence on which it has a positive limit. For $\varepsilon>0$, we use Sobolev's inequality and the first identity in~\eqref{scalardefect} to show that the residual action is at least
\[
b_\varepsilon:=\frac1N S^{\frac{N}{2}}\varepsilon^{-\frac{N-2}{2}}.
\]
Require also
\begin{equation}\label{epsbubble}
\varepsilon\leq\left(\frac{S^{\frac{N}{2}}}{2NL}\right)^{\frac{2}{N-2}}.
\end{equation}
Then $b_\varepsilon\geq2L$, whereas $c_\varepsilon\leq L-1$. Using the second identity in~\eqref{scalardefect}, we rule out that residual. For $\varepsilon=0$, we obtain $R_n\to0$ directly from the first identity. Since $\omega>0$ is fixed,
\[
\|W_n-W\|_{H^1}^2
\leq\max\{1,\omega^{-1}\}R_n\longrightarrow0.
\]
Thus $W_n$ converges strongly in $H^1$. By~\eqref{scalar-nonzerolimit}, its limit changes sign and attains the minimum. By the deformation argument, every minimizer is critical, and by radial averaging we recover the full scalar equation. Using the same argument, we also obtain compactness when the frequencies converge in $K_r$.

\emph{A uniform bound for the change in the action level.}
For an attained minimizer, the projection onto the pure-$r$ constraints has amplitudes
\[
t_\pm=(\frac{A_\pm}{B_\pm})^{\frac{1}{r-2}}
=\left(1-\frac{\varepsilon \|W^\pm\|_s^s}{A_\pm}\right)^{-\frac{1}{r-2}}.
\]
Since $0\leq\varepsilon \frac{\|W^\pm\|_s^s}{A_\pm}\leq\frac{1}{2}$, we differentiate $(1-x)^{-\frac{2}{r-2}}$ on $[0,\frac{1}{2}]$ and apply the mean value theorem to obtain
\[
0\leq t_\pm^2-1
\leq\frac{2}{r-2}2^{\frac{r}{r-2}}
\frac{\varepsilon \|W^\pm\|_s^s}{A_\pm}.
\]
The projected pure action is $k_r\sum_\pm A_\pm t_\pm^2$, while
$J_{\omega,\varepsilon}(W)=k_r\sum_\pm A_\pm+(k_s-k_r)\varepsilon\sum_\pm \|W^\pm\|_s^s$.
Consequently the positive constant
\[
C_0:=\frac{2^{\frac{r}{r-2}}}{r}S^{-\frac{s}{2}}M^{\frac{s}{2}}
\]
satisfies
\begin{equation}\label{scalarOeps}
0\leq c_0(\omega)-c_\varepsilon(\omega)
\leq k_r\sum_\pm A_\pm(t_\pm^2-1)
\leq C_0\varepsilon.
\end{equation}
In the last step we used $\sum_\pm \|W^\pm\|_s^s\leq S^{-\frac{s}{2}}M^{\frac{s}{2}}$.

\emph{Choosing the frequency and imposing mass one.}
Define $F_\varepsilon(\omega)=c_\varepsilon(\omega)-\frac{\omega}{2}$ for $\varepsilon\geq0$. In particular $F_0''(\omega)=\sigma(\sigma+1)d_r\omega^{\sigma-1}>0$, and $F_0'(\omega_0)=0$. The boundary gap
\[
\Delta_r:=\min\{F_0(\frac{\omega_0}{2}),F_0(2\omega_0)\}-F_0(\omega_0)
\]
is therefore positive. Choose
\begin{equation}\label{epschoice}
\varepsilon_*=
\min\left\{1,\frac{S^{\frac{s}{2}}}{2M^{\frac{s-2}{2}}},
\left(\frac{S^{\frac{N}{2}}}{2NL}\right)^{\frac{2}{N-2}},
\frac{\Delta_r}{2C_0}\right\}>0.
\end{equation}
All quantities in this choice depend only on $N,r$.

To compare nearby frequencies, fix a minimizer $W$ at $\omega$. At $\omega+h$, its two projection amplitudes solve
\[
A_\pm+h\|W^\pm\|_2^2
=B_\pm t_\pm(h)^{r-2}
+\varepsilon \|W^\pm\|_s^s t_\pm(h)^{s-2}.
\]
The derivative of the right side at $t=1$ is
$(r-2)B_\pm+(s-2)\varepsilon \|W^\pm\|_s^s>0$, so by the implicit-function theorem we have $t_\pm(h)=1+O(h)$. The amplitude derivatives of the action vanish at $h=0$. Expanding the projected action, we consequently obtain
\begin{equation}\label{scalarenv}
c_\varepsilon(\omega+h)\leq c_\varepsilon(\omega)
+\frac h2\|W\|_2^2+O(h^2).
\end{equation}
From the estimate, we obtain upper semicontinuity. Using compactness of minimizers, we obtain lower semicontinuity and conclude that $c_\varepsilon$ is continuous on $K_r$.

For $0\leq\varepsilon<\varepsilon_*$ and $\omega\in\partial K_r$, we combine~\eqref{scalarupper} and~\eqref{scalarOeps} to obtain
\[
F_\varepsilon(\omega)
\geq F_0(\omega_0)+\Delta_r-C_0\varepsilon
>F_\varepsilon(\omega_0).
\]
Thus $F_\varepsilon$ has an interior minimum. At such a minimum, we obtain from~\eqref{scalarenv}, for any ground state $W$,
\[
0\leq F_\varepsilon(\omega+h)-F_\varepsilon(\omega)
\leq\frac h2(\|W\|_2^2-1)+O(h^2).
\]
For $h=t>0$ divide by $t$; for $h=-t$ again divide by the positive number $t$. Letting $t\downarrow0$ yields both $\|W\|_2^2-1\geq0$ and $\|W\|_2^2-1\leq0$. Hence the mass is exactly one. Its energy equals $F_\varepsilon(\omega)$ and satisfies
\[
F_\varepsilon(\omega)\leq F_0(\omega_0)
=-\frac{\sigma\omega_0}{2(\sigma+1)}<0.
\]
\end{proof}

\begin{theorem}\label{thm:sobpure}
Assume~\eqref{reverse-assumptions} with $s=2^*$ and the balance~\eqref{balance}. With $\delta_r$ from~\eqref{deltar}, define
\begin{equation}\label{Ddef}
C_{N,r}=\varepsilon_*(N,r)^{\frac{2-\delta_r}{s-r}}>0,
\end{equation}
where~\eqref{epschoice} is an admissible choice of $\varepsilon_*$. If the scalar target mass $m$ satisfies
\begin{equation}\label{massinterval}
0<m<C_{N,r}\bar\mu^{-\frac{2-\delta_r}{s-r}}
\kappa^{-\frac{s-2}{s-r}},
\end{equation}
there exists $\lambda>0$ and a solution in $\mathcal Y_{\alpha\beta}$ with masses $(a,b)=(\alpha m,\beta m)$ and both components sign-changing. At $(\lambda,\lambda)$ it minimizes the action on the entire intersection of the four-sign Nehari set with $\mathcal Y_{\alpha\beta}$, and hence among all radial nodal solutions in that proportionality class. Its energy is negative. Every fixed positive mass pair on this ray is therefore covered for sufficiently small $\nu>0$.
\end{theorem}
\begin{proof}
For a target mass $m>0$, choose the positive spatial scale
\begin{equation}\label{criticalscaling}
\rho=(\kappa m^{\frac{r-2}{2}})^{\frac{1}{2-\delta_r}},\qquad
w(x)=\sqrt m\,\rho^{\frac{N}{2}}W(\rho x),\qquad
\lambda=\rho^2\omega.
\end{equation}
With the change of variables $y=\rho x$, we obtain
\[
\|w\|_2^2=m\|W\|_2^2,\qquad
\|\nabla w\|_2^2=m\rho^2\|\nabla W\|_2^2.
\]
We divide the equation for $w$ by $\sqrt m\,\rho^{\frac{N}{2}+2}$. By the definition of $\rho$,
the lower-power coefficient is $\kappa m^{\frac{r-2}{2}}\rho^{\delta_r-2}=1$.
For the critical coefficient, $\frac{N(s-2)}{2}-2=s-2$ because $s=2^*$. Thus
\begin{align}
\varepsilon
&=\bar\mu m^{\frac{s-2}{2}}\rho^{s-2}\nonumber\\
&=\bar\mu\kappa^{\frac{s-2}{2-\delta_r}}
 m^{\frac{s-2}{2}+\frac{(r-2)(s-2)}{2(2-\delta_r)}}\nonumber\\
&=\bar\mu\kappa^{\frac{s-2}{2-\delta_r}}m^{\frac{s-r}{2-\delta_r}}.
\label{epsilonformula}
\end{align}
For the last equality, the mass exponent is
$\frac{(s-2)(r-\delta_r)}{2(2-\delta_r)}$, and
$(s-2)(r-\delta_r)=2(s-r)$ when $s=\frac{2N}{N-2}$.
Since $2-\delta_r>0$ and $s-r>0$, we recover exactly~\eqref{massinterval} by solving $\varepsilon<\varepsilon_*$ for $m$.

By Lemma~\ref{lem:scalarunit}, we select $\omega>0$ and a radial nodal $W$ of mass one. The rescaled function $w$ has mass $m$ and both of its sign parts are nonzero. Hence $(\sqrt\alpha\,w,\sqrt\beta\,w)$ has masses $(\alpha m,\beta m)$ and both components change sign.

With the same changes of variables, we find
\[
I_{\lambda,\bar\mu,\kappa}(w)
=m\rho^2J_{\omega,\varepsilon}(W).
\]
With this scaling, we obtain a bijection on the radial space that preserves each nonzero sign and maps the two nodal constraints bijectively. The action is multiplied by a positive factor, so we preserve the entire nodal minimum. Thus $w$ has mass $m$ and the scalar least-action property. Using the reduction~\eqref{actionreduction}, we then obtain the full system, its two masses, and the asserted minimum in $\mathcal Y_{\alpha\beta}$. Subtracting half the frequency times the mass in the action identity, we obtain the same positive multiplier for the energy and conclude that the energy is negative.
\end{proof}

\subsection{Mass critical powers: a different interval}
Now take $s=p_c=2+\frac{4}{N}$ and retain $2<r<p_c$. Define
\begin{equation}\label{nodalmasscritical}
M_N^{\nod}=2d_{p_c,1}(1).
\end{equation}
By~\eqref{masslaw}, this is the squared $L^2$ norm of every  scalar radial nodal action ground state with coefficient one. It differs from the mass of a positive ground state.

\begin{theorem}\label{thm:masspure}
Assume~\eqref{reverse-assumptions} with $s=p_c$ and~\eqref{balance}. For every $\nu>0$ and
\begin{equation}\label{masscriticalinterval}
0<m<\bar\mu^{-\frac{N}{2}}M_N^{\nod},\qquad(a,b)=(\alpha m,\beta m),
\end{equation}
there is a normalized solution with positive equal frequencies and both components sign-changing. It is an action ground state within $\mathcal Y_{\alpha\beta}$ at its frequencies, and its energy is negative.
\end{theorem}
\begin{proof}
\emph{Attainment and variation of the scalar level.}
For each $\lambda>0$, let $c(\lambda)$ be the infimum of the scalar action~\eqref{reducedaction} on the radial two-sign constraints~\eqref{scalarmixedN}. Both powers satisfy $2<r<s<2^*$. For each sign, the corresponding amplitude equation has a unique positive root because its nonlinear side is strictly increasing. On the constraints,
\[
I_{\lambda,\bar\mu,\kappa}(w)
=k_s\bar\mu\|w\|_s^s+k_r\kappa\|w\|_r^r
\geq k_r\|w\|_{1,\lambda}^2.
\]
Using a fixed nodal trial, we obtain a finite upper bound. For each sign, put $A=\|w^\pm\|_{1,\lambda}^2$. From the constraint and the two embeddings, we obtain
$A\leq C(A^{\frac{s}{2}}+A^{\frac{r}{2}})$. Dividing by $A>0$ rules out $A\to0$. On the nodal minimizing sequence, we then use the constraint to bound the sum of the two sign nonlinearities away from zero. By a two-amplitude deformation, we construct a bounded Palais--Smale sequence $w_n$ close to that minimizing sequence in $H^1$. Truncation is Lipschitz in $L^r,L^s$, so for each sign the lower bound passes to this sequence:
\[
\bar\mu\|(w_n)^\pm\|_s^s+\kappa\|(w_n)^\pm\|_r^r\geq c>0.
\]
After extracting a weakly convergent subsequence, we use radial compactness at both powers to obtain strong $L^r,L^s$ convergence to a function $w$. The displayed sum converges to
$\bar\mu\|w^\pm\|_s^s+\kappa\|w^\pm\|_r^r\geq c$.
Therefore both signs of the weak limit are nonzero. We pass to the weak equation and test its difference from the Palais--Smale equations to obtain
$\|w_n-w\|_{1,\lambda}^2\to0$. Since $\lambda$ is fixed and positive, we obtain strong $H^1$ convergence and attainment, including strong $L^2$ convergence. Every minimizer is a full scalar solution by the same deformation and radial averaging arguments.

These bounds are uniform on compact positive frequency intervals. Thus minimizers are compact along convergent frequencies. Using the two-amplitude projection at $\lambda+h$, with $h\in\R$ small, we obtain, for any minimizer $w$ at $\lambda$,
\begin{equation}\label{masspureenv}
c(\lambda+h)\leq c(\lambda)+\frac h2\|w\|_2^2+O(h^2).
\end{equation}
As in~\eqref{scalarenv}, we solve the two scalar amplitude equations and expand at their stationary values. Using compactness as well, we obtain continuity of $c$.

\emph{The two frequency endpoints.}
Use a pure-$r$ nodal ground state as a trial and take the maximum along its two sign amplitudes. Removing the nonpositive pure-$s$ term can only increase that maximum. By scalar dilation, we therefore obtain
\[
0<c(\lambda)\leq d_{r,\kappa}(1)\lambda^{1+\frac{2}{r-2}-\frac{N}{2}}.
\]
Since $\frac{2}{r-2}-\frac{N}{2}>0$, it follows that $\frac{c(\lambda)}{\lambda}\to0$ as $\lambda\downarrow0$.

For the other endpoint, introduce a rescaled function $Z$ by $w(x)=\lambda^{\frac{N}{4}}Z(\sqrt\lambda x)$, which preserves the mass because $s=p_c$. Set
$\eta=\kappa\lambda^{-\frac{2-\delta_r}{2}}$, which tends to zero as $\lambda\to\infty$. Dividing the action by $\lambda$ gives
\[
\lambda^{-1}I_{\lambda,\bar\mu,\kappa}(w)
=\frac12\|Z\|_{1,1}^2-\frac{\bar\mu}{s}\|Z\|_s^s
-\frac\eta r\|Z\|_r^r.
\]
Here $\|Z\|_{1,1}^2=\|\nabla Z\|_2^2+\|Z\|_2^2$. With this scaling, we map the nodal constraint sets bijectively. Hence $\frac{c(\lambda)}{\lambda}$ is the nodal minimum of the displayed expression. Using a pure-$s$ ground state, we obtain $\frac{c(\lambda)}{\lambda}\leq d_{s,\bar\mu}(1)$. By coercivity with coefficient $k_r$, we bound the scaled minimizers $Z$ uniformly in $H^1$.

For such a minimizer, write $A_\pm=\|Z^\pm\|_{1,1}^2$ and $B_\pm=\bar\mu\|Z^\pm\|_s^s$. Its constraint is $A_\pm=B_\pm+\eta\|Z^\pm\|_r^r$. The uniform norm bound and $r>2$ give
$\eta\|Z^\pm\|_r^r\leq C\eta A_\pm$, so $B_\pm\geq \frac{A_\pm}{2}$ for small $\eta$. The scalar lower-bound argument then bounds both $A_\pm,B_\pm$ away from zero. Projection onto the pure-$s$ constraints uses amplitudes
$t_\pm=(\frac{A_\pm}{B_\pm})^{\frac{1}{s-2}}=1+O(\eta)$.
All coefficients are bounded, so the projected pure action differs from the scaled mixed action by $O(\eta)$. Hence
\[
d_{s,\bar\mu}(1)\leq\frac{c(\lambda)}\lambda+C\eta
\leq d_{s,\bar\mu}(1)+C\eta.
\]
Finally, multiplication of a power solution by $\bar\mu^{-\frac{1}{s-2}}$ multiplies its action by $\bar\mu^{-\frac{2}{s-2}}=\bar\mu^{-\frac{N}{2}}$. We have proved
\begin{equation}\label{masspureasymptotic}
\lim_{\lambda\to\infty}\frac{c(\lambda)}\lambda
=\bar\mu^{-\frac{N}{2}}d_{p_c,1}(1).
\end{equation}

\emph{Exact mass selection.}
For $m$ in~\eqref{masscriticalinterval}, set $F(\lambda)=c(\lambda)-\frac{m\lambda}{2}$. The function $F$ is continuous, tends to zero as $\lambda\rightarrow0^+$, and is negative for small positive $\lambda$. By~\eqref{masspureasymptotic}, it tends to $+\infty$ as $\lambda\to\infty$. Its global minimum is therefore attained at a positive frequency.

For every ground state $w$ at that frequency, we obtain from~\eqref{masspureenv}
$0\leq F(\lambda+h)-F(\lambda)\leq \frac{h}{2}(\|w\|_2^2-m)+O(h^2)$.
Taking $h=t$ and $h=-t$, dividing by $t>0$ and letting $t\downarrow0$, we obtain $\|w\|_2^2=m$. Its energy is the negative minimum of $F$. Using the exact reduction and~\eqref{actionreduction}, we obtain the vector conclusion.
\end{proof}

\section*{Authors' addresses}

\noindent\textbf{Yumeng Chou}\\
School of Mathematics and Statistics, Nantong University, Nantong, Jiangsu, China\\
Email: \texttt{chouyumeng@stmail.ntu.edu.cn}

\medskip

\noindent\textbf{Huabo Zhang}\\
School of Mathematics and Statistics, Nantong University, Nantong, Jiangsu, China\\
Email: \texttt{lutzhb@163.com}

\end{document}